\documentclass[12pt]{amsart} 
\newfont{\sheaf}{eusm10 scaled\magstep1}

\newcommand{\C}{\ensuremath{\mathbb{C}}}
\newcommand{\D}{\ensuremath{\mathbb{D}}}
\newcommand{\R}{\ensuremath{\mathbb{R}}} 
\newcommand{\Z}{\ensuremath{\mathbb{Z}}}
\newcommand{\F}{\ensuremath{\mathbb{F}}} 
\newcommand{\Q}{\ensuremath{\mathbb{Q}}}
\newcommand{\N}{\ensuremath{\mathbb{N}}} 
\newcommand{\hol}{\ensuremath{\mathcal{O}}}
\newcommand{\HH}{\ensuremath{\mathbb{H}}} 
\newcommand{\PP}{\ensuremath{\mathbb{P}}} 
\newcommand{\RR}{\ensuremath{\mathcal{R}}}
\newcommand{\A}{\ensuremath{\mathcal{A}}}

\newcommand{\SSS}{\ensuremath{\mathcal{S}}}
\newcommand{\FF}{\ensuremath{\mathcal{F}}}
\newcommand{\HHH}{\ensuremath{\mathcal{H}}}
\newcommand{\PPP}{\ensuremath{\mathcal{P}}}

\newcommand{\LL}{\ensuremath{\mathcal{L}}}

\newcommand{\Proof}{{\it Proof. }} 
\newcommand{\ra}{\ensuremath{\rightarrow}}

\newtheorem{teo}{Theorem}[section] 
\newtheorem{df}[teo]{Definition} 
 
\newtheorem{cor}[teo]{Corollary} 
\newtheorem{oss}[teo]{Remark} 

\newtheorem{ex}[teo]{Example}

\newtheorem{problem}[teo]{Problem}

\def\eea{\end{eqnarray*}}
\def\bea{\begin{eqnarray*}}

\begin{document}

\title{From Abel's Heritage: transcendental objects in algebraic geometry
and their algebrization.}

\author{Fabrizio Catanese\\
 Universit\"at Bayreuth
 }

\footnote{
The present research, 
an attempt  to treat history and
sociology of mathematics and  mathematics all at the same time, took place
in the framework of the Schwerpunkt "Globale Methode in der komplexen Geometrie",
and of the EC Project EAGER.   }
\date{October 29, 2002}
\maketitle

{\bf Welcome the Abel Prize and  long live the memory of Abel, 

long live mathematics! \footnote {Was the greatest enemy of mathematics
Alexander the Great, who cut the Gordian knot instead of peacefully writing a
book  about it? }}\\

\section{Introduction}

$${\bf CONTENTS}$$
\begin{itemize}

\item
2. Abel, the algebraist ?

\item
3. The {\bf geometrization} of Abel's methods.
\item
4. Algebraization of the Geometry.
\item
5. Further\footnote{ This topic has in fact been treated quite extensively in the
contribution by Ciliberto.} links to the Italian School.
\item
6. More new results and open problems.
\end{itemize}

The Encyclopedic Dictionary of Mathematics (edited in Japan) is not
a book devoted to the history of mathematics, but  tries instead to
briefly introduce the reader to the current topics of mathematical
research.  By non only lexicographical coincidence it starts with 
"Abels, Niels Henrik" as topic 1.

It contains a succint biography of Abel

" Niels Henrik Abel (august 5, 1802 - april 6 , 1829) ... 
In 1822, he entered the University of Christiania [today's Oslo]
.. died at twenty-six of tuberculosis.
His best known works are: the result that algebraic equations of order
five or above cannot in general be solved algebraically; the result
that \footnote{ A star next to a theme denotes that a section of the
dictionary is devoted to the discussion of the topic}*Abelian equations
[i.e., with Abelian Galois group]
can be solved algebraically; the theory of *binomial series and of
*elliptic functions; and the introduction of *Abelian functions. His
work in both algebra and analysis, written in a style conducive to easy
comprehension, reached the highest level of attainment of his time."

 Talking  about Abel's heritage entails thus
 talking about a great part of modern mathematics, as it is shown by the 
 ubiquity of concepts such as Abelian  groups, Abelian  integrals
 and functions, Abelian Varieties, and of relatives of theirs as anabelian
 geometry, nonabelian Hodge theory...
 
 Writing is certainly a more difficult task than talking, when the time limits
force us to plan our way on secure and direct tracks: for this reason we
decided that the present text, with the exception of a couple of  protracted
mathematical discussions, should essentially be the text of our oral exposition
at the Abel Bicentennial Conference. Thus its aim is just to lead the reader along
quite personal views on history and development of mathematics, and on certain topics
in the still very alive subject of transcendental algebraic geometry.

One could declare its Leitmotiv to be  G.B. Vico's theory of cycles in the
history of mankind, adapted to the analysis of mathematical evolutions and
revolutions:

Geometry in ancient Greece, Algebra by the Arabs and in early
Renaissance, Geometry again by B. Cavalieri and his indivisibles,
Analysis and Physics by the Bernoulli's,.. and later on an intricate
 succession of  points of view and methods, often alternative to each
other, or striving in directions opposite one to the other, which  all
together enriched our knowledge and understanding of the mathematical reality.

Therefore, if we conceive algebra, geometry ..  more as methodologies than
as domains of knowledge, comes out naturally the difficult question: which
 way of doing mathematics is the one we are considering ?
 
 This question, probably a sterile question when considering the history of
mathematics, is however a very important one when we are making choices for
future directions of mathematical research:  to purport this assertion it will
suffice only to cite the (for me, even exaggerated) enthusiasm of nowadays
algebraic geometers for the new insights coming to their field by physical
theories, concepts and problems.

In any case, in our formerly bourgeois world, idle questions with
provocative answers to be defended at tea time at home or in a Cafe', have
often motivated interesting discussions, and it is just my hope to be
able to do the same thing here.

\section{Abel, the algebraist?}

So I will start, as due, by  citing Hermann Weyl's point of view
(cf. \cite{Weyl31},\cite{Weyl} and also \cite{yag}, pp. 26 and 151, for comments),
expressed in an address directed towards mathematics teachers, and later
published in the Journal : Unterrichtsbl\"atter f\"ur Mathematik und
Naturwissenschaften, Band 38 (1933), s. 177-188.

There are two Classes of mathematicians:

\begin{itemize}
\item
ALGEBRAISTS: as Leibniz, Weierstrass
\item
GEOMETERS-PHYSICISTS: as Newton, Riemann, Klein

\end{itemize}

and people belonging to different classes may tend to be in conflict with each
other.

The tools of the algebraists are : logical argumentation, formulae and
their clever manipulation, algorithms.

The other class relies more on intuition, and graphical and visual
impressions. For them it is more important to find a new truth than an
elegant new proof.

The concept of rigour is the battlefield where the opposite parties 
confront themselves, and the conflicts which  hence derived were 
sometimes harsh and longlasting.

The first well known example is the priority conflict between Leibniz and
Newton concerning the invention of the Calculus (which however was invented
independently by the two scientists, as it is currently agreed upon).

The  inputs which the two scientists provided did indeed integrate themselves
perfectly. On one hand the pure algebraic differential quotient  $dy/dx$ would
be a very dry  concept (algebraists nevertheless are still nowadays very keen
on inflicting on us the abstract theory of derivations!) without the intuition
of velocities and curve tangents; on the other hand, in the analysis of
several phenomena, a physical interpretation of Leibniz's rule  can turn out to
be amazingly complicated .

More closely related with Abelian integrals and their periods was 
Weierstrass' constructive criticism of the "Principles" by Riemann and
Dirichlet.

As also pointed out in the contribution\footnote{ Here and after, I will refer to
the oral contributions given at the Abel Bicentennial Conference, and not to the
articles published in this  Volume. } by Schappacher, this conflict soon
became the Berlin-G\"ottingen conflict, and almost deflagrated between Weierstrass and
Felix Klein (who continued on the way started as a student of Clebsch)\footnote{   We
have noticed that another article devoted to this topic has appeared after we gave the
talk,  namely
\cite{bott02} by U. Bottazzini.} .

Klein's antipathy for Weierstrass was more intellectual than personal:
Klein put a special enphasys on geometrical and physical intuition,
which he managed to develop in the students by letting them construct
solid (plaster,   or metal) models of curves and surfaces in
ordinary 3-space, or letting them draw very broad (1 x 2 meters) paper
tables of cubic plane curves with an explicit plotting of their
$\Q$-rational points.

A concrete witness to this tradition is the exhibition of plaster
models of surfaces which are to be found still nowadays in the Halls
of the mathematical Institute in G\"ottingen. These models were then
produced by the publishing company L. Brill in Darmstadt, later by the 
Schilling publishing company, and sold
around the world: I have  personally seen many of those 
in most of the older Deparments I have visited (cf. the 2 Volumes edited by
G.Fischer on  "Mathematical Models" \cite{fisch86}
\footnote{ Writes G. Fischer : " There were certainly other reasons than
economic for the waning interest in models. ...... More and more general
and abstract viewpoints  came to the forefront of mathematics....... Finally
Nicolas Bourbaki totally banned pictures from his books."} ).

A similar trick, with computer experiments replacing the 
construction of models, is still very much applied nowadays in the case where
professors have to supervise too many more students' theses than they can
really handle.

Weierstrass had a victory, in the sense that not only the theory of
calculus, but also the theory of elliptic functions is still nowadays
taught almost in the same way as it was done in his Berlin lectures.

But Klein's "defeat" (made harder by the long term
competition with Poincare' about the proof of the uniformization theorem, see later)
was however a very fertile humus for the later big growth of the G\"ottingen
influence, and certainly Weyl's meditations which we mentioned above were reflecting
this  historically important controversy.

Where does then Abel stay in this classification?
I already took position, with my choice of the title of this section: Abel is for
me an algebraist and I was glad 
\footnote {For many years on, until after the first world war, there
used to be a course in Italian Universities entitled "Lezioni di analisi
algebrica ed infinitesimale". This shows that  the birth of
Analysis as a new separated branch, trying to appear more on the side of
applied mathematics, is a relative novelty which ends the reconciliation
made for the dualism Leibniz-Newton.  } 
to hear  Christian Houzel stressing in his
contribution  the role of Abel's high sense of rigour. Abel's
articles on the binomial coefficients,  on the summation of series, and
on the solution of algebraic equations testify his deep concern for
the need of  clear and satisfactory proofs
\footnote { cf. the article by: J J O'Connor and E F Robertson in
http://www-history.mcs.st-andrews.ac.uk/References/Abel.html, citing his letter
to Holmboe from Berlin. Here we can read : " In other words, the most
important parts of mathematics stand without foundation".}.

Of course, like many colours
are really a  mixture of pure colours, the same occurs for mathematicians,
and by saying that he was deep down an "algebraist" I do not mean to
deny that he possessed a solid geometrical intuition, as we shall later
point out.

In fact, Abel himself was proud to introduce himself  during his travels as
'Professor of Geometry' \footnote{To be perfectly honest, all we know is that in 1826
he signed himself in at the ``Goldenes Schiff'' in Predazzo as `Abel, professore della
geometria.'}.

The best illustration of his synthetical point of view is shown by the
words (here translated from French) with which he begins his Memory (XII-2
"M\'emoire sur une propri\'et\'e g\'en\'erale d'une classe tres \'etendue de
fonctions transcendantes."):
\footnote{  This important article was lying at the centre of the contribution
by Griffiths, and is also amply commented upon in Kleiman's contribution.}

" The transcendental functions considered until nowadays by the
geometers are a very small number. Almost all the theory of
transcendental functions is reduced to the one of logarithmic, exponential and
circular functions, which are essentially one only kind. Only in recent times one
has begun to consider some other functions. Among those, the elliptic
transcendentals, about which M. Legendre developed so many elegant and
remarkable properties, stay in the first rank."

The statement "are essentially one only kind" is the one we want now to comment
upon.

Algebraists like indeed short formulae, and these are in this case
available. It suffices to consider the single formula : $exp (x) :=
\Sigma_{n=0}^{ \infty} \frac{ x^n} {n!} $, and then, by considering $ exp(ix)$
and the  inverse functions of the ones we can construct by easy algebraic
manipulations, we get easily ahold of the wild proliferation of functions
 which for instance occupy the stage of the U.S.A. Calculus courses 
 (sin, cos, sec, cosec, tg, cotg, and their hyperbolic analogues).

No doubt, synthesis is a peculiarity of pure mathematics, and applied or
taught mathematics may perhaps need so many different names and functions: but,
 pretty sure,  Abel stood by the side of synthesis and
conciseness.

His well known saying, that he was able to learn so rapidly because he had been: "
studying the Masters and not their pupils", is quite valid nowadays.
Today there is certainly 
an inflation of books and divulgations, many are second,  third  hand  or even
further. Abel's point of view should be seriously considered by some pedagogists
who want to strictly regulate children's learning, forcing them to study n-th hand
knowledge. Perhaps this is a strictly democratic principle, by which one wants to
prevent some children from becoming precociously wise (as Abel did), and possibly
try to stop their intellectual growth (this might be part of a more general ambitious
program, sponsored by Television Networks owners).

In any case, Abel had read the masters, and he knew many functions: he still
belongs to the mathematical era where functions are just concretely defined
objects and not subsets of a Cartesian product satisfying a geometrical
condition. One of Abel's main contributions was to consider not only ample classes
of functions defined by integrals of algebraic functions, but to study their
inverse functions and their periodicity (the so called elliptic functions being
among the latter).

At his time did not yet exist the concept of "Willk\"urliche Funktionen"
(= "arbitrary functions"),
quite central e.g. for Weierstrass, Dini, Peano  and Hilbert (cf. \cite{dini},
\cite{weier}, \cite{g-p} )
 and which motivated much of the developments in the theory of sets leading
to the construction of several pathological situations (as Lebesgue's non constant
function with derivative almost everywhere zero, \cite{leb} ).

Most of the functions  he considered were in fact written as $\int_{x_0}^x y(t)
dt $, where $y(x)$ is an {\bf algebraic function} of $x$, which simply
means that the function is defined on some interval in $\R$ and that there
exists a polynomial
$ P(x,y) \in \C[x,y]$ \footnote{ Observe that we wrote $\C[x,y]$ instead of
$\R[x,y]$. It is commonly agreed that, if we want to summarize in two words
which the greatest contribution of Abel and Jacobi was, then it was to consider the
elliptic integrals not just as functions of a real variable, but also as functions of a
complex variable. So, we owe to them the birth of the theory of holomorphic functions.}
such that 
$ P(x,y(x)) \equiv 0.$ The functions given by these integrals, or by sums of
several of these, are nowadays called {\bf Abelian functions}.

The above statement  is by and large true, with however a single important exception,
 concerning Abel's treatment of functional equations: there he considers quite
generally the functions which occur as solution of certain {\bf functional
equations}.

As example, we take the content of an article also considered in Houzel's talk,
VI-Crelle I (1826). This article is a real gem: it anticipates S.Lie's
treatment of Lie group germs, and yields actually a stronger result (although, under
the assumption of commutativity):
\begin{teo}
Let $f : U \ra \R$ be a germ of function defined on a neighbourhood of
the origin, $ \R^2 \supset U \ni 0 $, such that
$$ f (z, f(x,y)) $$
is symmetric in $x,z,y$ (i.e., in today's terminology, we have an Abelian
Lie group germ in 1-variable). Then there exists a germ of change of variable $
\psi: (\R,0) \ra (\R,0) $ such that our Lie group becomes $(\R, +)$, or, more
concretely, such that 
$$ \psi (f(x,y)) = \psi(x) + \psi(y). $$

\end{teo}

Among the Masters' work which Abel studied was certainly, as already
mentioned, M. Legendre and his theory of elliptic integrals. For these,
already considered by Euler and Lagrange, 
Legendre devised a  {\bf normal form} (here $R(x)$ is a rational
function of $x$): 

$$ \int_{x_0}^t \frac{R(x)} {\surd \overline{ (1 - x^2) ( 1- k^2 x^2) }
} dx.$$

In the remarkable paper XVI-1 (published on Crelle, Bd. 2,3 (1827,
1828)), entitled "Recherches sur les fonctions elliptiques", Abel, as we
already mentioned , writes clearly, after observing that the study of
these elliptic integrals can be reduced to the study of integrals of
the first, second and third kind \footnote{Since
an elliptic integral is the one where we consider a square root  
$\surd\overline{ P(x) } $ where $P$ is a polynomial of degree $3$ or
$4$, one can reduce it, after applying a projective transformation of
the line $\PP^1$,  to a square root $\surd \overline{ Q(y^2) }$, where
$Q$ is quadratic, and then we can view it as an integral on the unit
circle, whose projection to the $x$-line  $\PP^1$ yields a double
cover.}
$$ \int \frac{ d \theta} {\surd \overline{ 1 - c^2 sin^2 \theta }  };
\int  d \theta \surd \overline{ 1 - c^2 sin^2 \theta }  , \int \frac{
d \theta } { ( 1 + n sin^2 \theta) \surd \overline{ 1 - c^2 sin^2 \theta
}  },$$

" These three functions are the ones that M. Legendre has considered,
especially the first, which enjoys the most remarkable and the simplest
properties. I am proposing myself, in this Memoir, to consider the
inverse function, i.e., the function
$\phi (a)$ , determined by the equations 
$$ a = \int \frac{ d \theta} {\surd \overline{ 1 - c^2 sin^2 \theta } 
}$$
$$ sin \theta = \phi (a) = x ."$$
In the following pages , where he manages to
give a simple proof of the double periodicity of the given function
$\phi$,  Abel shows the clarity  of his geometric intuition. 

He simply observes that in Legendre's normal form one should win the
natural resistance to consider non real roots, and actually it is much
better to consider the case where $ k^2 < 0$, and thus he considers
(his notation) the integral 
$$ \int_{x_0}^t \frac{1} {\surd \overline{ (1 - c^2 x^2) ( 1+ e^2 x^2) }
} dx,$$ where $c,e$ are strictly positive real numbers.

The roots of the radical are the points $ \pm 1/c , \pm (1/e) \surd
\overline{-1}$, we have a rectangular symmetry around the origin
and we have two periods $  \omega , \tilde{ \omega }$ obtained by integrating on
the two closed paths lying over the segments joining pairs of opposite roots:
$$\omega = 4 \int_{0}^{1/c} \frac{1} {\surd \overline{ (1 - c^2 x^2) (
1+ e^2 x^2) } } dx ,\ \tilde{ \omega } = 4 \int_{0}^{i/e} \frac{1}
{\surd
\overline{ (1 - c^2 x^2) ( 1+ e^2 x^2) } } dx.$$
It is straightforward to observe that the two periods (of the {\bf real
curve}) 
are such that $\omega \in \R$,  respectively $\tilde{ \omega } \in i \R$,
thus they are linearly independent over the real numbers.

The conclusion is that the inverse function is {\bf doubly periodic}
with one real and one imaginary period, so that its fundamental domain
is a rectangle with sides parallel to the real, resp. imaginary axis.
Since every elliptic integral of the first kind can be reduced to this
form, their inverse functions are all doubly periodic (unlike the
circular functions, which possess only one period).

Abel does not bother to highlight the geometry underlining his argument,
it is clear that he has developed a good geometrical intuition,
but his style is extremely terse and concise. This conciseness becomes
almost abrupt in the other article  XIII-Vol. 2 \footnote{ The second
volume of the edition \cite{abel81} by Sylow and Lie, Christiania 1881,
contains the
 unpublished papers of Abel, with a few exceptions. As the editors
remark, this edition, posterior by more than 30 years to the edition
of 1839 edited by the friend and colleague of Abel, Holmboe, was
financed by the Norwegian Parliament after the great demand
for Abel's works (Holmboe's edition went rapidly out of print) of which
especially the French Mathematical Society made itself interpreter. The
two editors decided to omit in the second volume three articles which
were partly based on an erroneous memoir of his youth, written in
Norwegian, where Abel thought he could prove that the general equation
of degree n can be solved by radicals. This is the only published
article which is  not appearing  in the Holmboe edition, nor in Volume I
of the edition by Sylow and Lie. }, "Theorie des transcendantes
elliptiques". 

This long memoir starts nulla interposita more (it was probably
unfinished): " For more simplicity I denote the radical by $\surd
\overline{R}$, whence we have to consider the integral $$\int \frac{P
dx} {\surd \overline{R} },$$ $P$ denoting a rational function of $x$."

It is divided into three chapters , the first devoted to the reduction
of elliptic integrals by means of algebraic functions,  the second to 
the reduction of elliptic integrals by means of logarithmic functions,
and finally the third is entitled " A remarkable relation which exists
among several integrals of the form 
$$\int \frac{dx} {\surd \overline{R} },\int \frac{x dx} {\surd
\overline{R} }, \int \frac{x^2 dx} {\surd \overline{R} }, \int \frac{dx}
{ (x-a)\surd \overline{R} }.$$

The memoir contains the explicit discussions of several concrete
problems concerning such reductions, and contains for many of those
problems explicit references to Legendre. It looks to me a rather early
work, because we directly see the influence of the study of Legendre,
but rather important for two reasons.

The first reason, already clear from the title of chapter III, is
that Abel here for the first time considers the question of the
relations holding among sums of elliptic integrals. This problem will be
considered more generally for all algebraic integrals of arbitrary
genus $g$ in his fundamental Memoir XII-1, entitled "M\'emoire sur une
propri\'et\'e g\'en\'erale d'une classe tres -\'etendue de
fonctions transcendantes", presented on october 10 1826 to the
Acad\'emie des Sciences de Paris, and published only in 1841.

The main theorem of the latter was so formulated: " If we have several
functions whose derivatives can be roots of the same algebraic
equation 
[if $y(x)$ is an algebraic function of $x$, i.e., there is a polynomial $P$
such that $P(x,y(x)) \equiv 0$, then for each rational function $f(x,y)$
there is a polynomial $F(x,y)$ such that $F(x, f(x,y(x))) \equiv 0$]\footnote{Here and
elsewhere, [..] stands for an addition of the present author}, with coefficients
rational functions of one variable [x], one can always express the sum of an arbitrary
number of such functions by means of an algebraic and a logarithmic function, provided
that one can establish among the variables of these functions a certain number of
algebraic relations".

The first theorem is then given through formula (12):
$$\int f(x_1,y_1) dx_1 + \int f(x_2,y_2) dx_2 + \dots \dots + \int f(x_{\mu},y_{\mu})
dx_{\mu} = v [(t_1, .. t_k)] \ :$$
here $f(x,y)$ is a rational function,  we take the $\mu$ points which form the
complete intersection of [$P(x,y) = 0$] and [$G_t(x,y) = 0$] where $G_t$ depends
rationally upon the parameter $ t = (t_1, .. t_k)$
, and the conclusion is, as we said,
that $v$ is the sum of a rational and of a logarithmic function.

Abel also explains clearly in the latter memoir that the number of these
relations is a number, which later on was called the genus of the curve
$C$ birational to the plane curve of equation
$P(x,y) = 0$. The way we understand the hypothesis of the theorem
nowadays is through the geometric condition: if the Abel sum of these
points is constant in the Jacobian variety of $C$. 
I will come back to the geometric interpretations in the next section,
let me now return to the second reason of importance of the cited
Memoir XIII, 2.

For instance, in Chapter I, Abel gives very explicit formulae, e.g. for
the reduction of integrals of the form 
$$\int \frac{x^m dx} {\surd \overline{R} },$$
where $R(x)$ is a polynomial of degree $3,4$, to the integrals
$$\int \frac{dx} {\surd \overline{R} },\int \frac{x dx} {\surd
\overline{R} }, \int \frac{x^2 dx} {\surd \overline{R} }.$$
The above integral is calculated by recursions, starting from the
equation 
$$  d (Q \surd\overline{R}) = S \frac{ dx} {\surd \overline{R} } ,$$ 
and writing explicitly
$$ S = \phi(0) + \phi(1) x + \dots  + \phi(m) x^m $$ 
$$ Q = f(0) + f(1) x + \dots  + f(m-3) x^{m-3} .$$
This is an example of Abel's mastery in the field of {\bf Differential
Algebra}. Although the modern reader, as well as Sylow and Lie, may
underscore the impact of these very direct calculations, it seems to me
that there has been a resurgence of this area of mathematics,
especially in connection with the development of computer
algorithms and programs which either provide an explicit integration of
a given function by elementary functions \footnote{I.e., by rational
functions or by logarithms (more generally, one can consider algebraic functions and
logarithms of 
these).}, or decide that the given function does not admit an integration
by elementary functions ( this problem was solved by Risch, and later concrete decision
procedures and algorithms were given by J Davenport and Trager  \cite{risch},
\cite{dav79-1},
\cite{dav79-2},\cite{dav81},
\cite{trag79}).

Differential algebra is also the main tool in the article (XVII-2)
" M\'emoire sur les fonctions transcendantes  de la forme $\int y \ dx$
ou $y$ est une fonction alg\'ebrique de $x$."

This paper looks very interesting and somehow gave me the impression
(or at least I liked to see  it in  this way) of being a forerunner
of the applications of Abelian integrals to 
questions of transcendence
theory.

This time I will state the main theorem by slightly altering Abel's
original notation

\begin{teo}
Assume that $\phi$ is a [non trivial] polynomial 
$$ \phi (r_1, \dots , r_{\mu}, w_1, \dots w_p) $$ and that $\phi \equiv
0$ if we set 
$$ r_i = \int y_i(x)\ dx, \ w_j = u_j (x), $$ where $y_i, u_j$ 
are algebraic functions.
Then there is a [non trivial] linear relation $ \Sigma_i  c_i\int y_i(x)\ dx
= P(x)$, with constant coefficients $c_i$ and with $P$ an algebraic function. 
\end{teo}

\begin{cor}
Let $y_i(x)\ dx$ be linearly independent Differentials of the I Kind on
a Riemann surface. Then the respective integrals are algebraically
independent.
\end{cor}

Nowadays, we would use the periodicity of these functions on the
universal cover of the algebraic curve, or the order of growth of the
volume of a periodic hypersurface (this argument was the one used
later by Cousin, cf. \cite{cous02}) to infer that the polynomial must be
linear (vanishing on the complex linear span of the periods). Abel's
argument is instead completely algebraic in nature, he chooses in fact
the polynomial
$\phi$ to be of minimal degree with respect to  $r_{\mu}$, and applies
$d / dx$ to the relation in order to obtain (Abel's
notation)\footnote{Why was not Abel using the notation $\partial /
\partial r_j$ ? }
$$ \Sigma_j \phi ' (r_j)\  y_j + \phi ' (x) = 0. $$
Abel shows then that $\phi$ has degree $1$ in $r_{\mu}$, since writing
$$ R = r_{\mu}^k + P_0  r_{\mu}^{k-1} + P_1  r_{\mu}^{k-2} +
\dots = 0$$ he gets:
$$ 0 = \frac{d}{dx} (R) =  r_{\mu}^{k-1} ( k  y_{\mu} + P_0') + \{ \dots \} 
r_{\mu}^{k-2} + \dots = 0 $$
whence, by the minimality of the degree $k$, this polynomial in $r_{\mu}$ has all
coefficients identically zero, in particular
$ ( k  y_{\mu} + P_0')
\equiv 0
$, therefore
$$r_{\mu} =
\int y_{\mu} dx = - (1/k) P_0$$ is the desired degree one relation.
By induction Abel derives the full statement that $\phi$ is linear in 
$r_1, \dots , r_{\mu}$.

\section{ The geometrization of Abel's methods.}

The process of geometrization and of a deeper understanding of Abel's
discoveries went a long way, with alternate phases, for over 150 years.
We believe that \footnote {As amply illustrated by Ciliberto in his
contribution.} a fundamental role for the geometrization was played by
the italian school of algebraic geometry, which then paved the way for
some of the more abstract developments in algebraic geometry.

Although it was very depressing for Abel that his fundamental Memoir
XII-1 was not read by A. Cauchy (this is the reason why it took more
than 15 years before it was published), still in Berlin Abel found the
enthusiastic support of L. Crelle, who launched his new Journal by
publishing the articles of Abel and Jacobi. Recall that finally fame
and recognition were reaching Abel through the offer of a professorship
in Berlin, which, crowning the joint efforts of Crelle and Jacobi,
 arrived however a few days after Abel had died.

Especially inspired by the papers of Abel  was Jacobi, in Berlin, who
also wrote a revolutionary article on elliptic function theory, entitled
"Fundamenta nova theoriae functionum ellipticarum". 

It was Jacobi who introduced the words 'Abelian integrals', 'Abelian
functions': Jacobi 's competing point of view (few 
competitions however
were so positive and constructive in the history of mathematics)
started soon to prevail.

Jacobi introduced the so called {\bf elliptic theta functions}, denoted
$\theta_{00}, \theta_{01},\theta_{10},\theta_{11},$ (cf. e.g. \cite{tric}
or \cite{mum3})
and expressed the elliptic functions, like Abel's $\phi$, inverse of the
elliptic integral of the first kind, as a ratio of theta functions.

A much more general definition of {\bf theta-series} (the expression
was later coined by Rosenhain and G\"opel, followers of Riemann) was 
given by B. Riemann, who defined his {\bf Riemann Theta function} as
the following series of exponentials
$$ \theta (z, \tau) : = \Sigma_{n \in \Z^g} \ exp ( 2 \pi \ i [
\frac{1}{2} \ ^t n \tau n + ^t n z ])$$
where $ z \in \C^g$, $\tau \in \mathcal H_g = \{ \tau \in Mat (g,g, \C) | \tau = ^t
\tau, Im (\tau) \ {\rm is \ positive \ definite} \}.$

The theta function converges because of the condition that $Im (\tau) 
$ is positive definite, it admits $\Z^g$ as group of periods,
being a Fourier series, and it has moreover a $\tau
\Z^g$-quasi-periodicity which turns out to be the clue for
constructing $2g$-periodic meromorphic functions as quotients of theta
series.

With a small variation (cf.
\cite{mum3}) one defines the theta-functions with characteristics 
$$ \theta [a, b](z, \tau) : = \Sigma_{n \in \Z^g} \ exp ( 2 \pi \ i [
\frac{1}{2} \ ^t (n+a) \tau (n+a) + ^t (n+a) (z+b) ]),$$
and the Jacobi functions $ \theta_{a, b}(z, \tau)$ are essentially
the functions $ \theta [a/2, b/2](2 z, \tau).$

Beyond the very explicit and beautiful formulae, what lies beyond this
apparently very analytic approach is the {\bf pioneeristic principle}
that any {\bf meromorphic function f} on a complex manifold $X$ can be
written as
$$ f = \frac{\sigma_1}{ \sigma_2} $$
of two relatively prime sections of a unique {\bf Line Bundle L } on
$X$.

This formulation came quite long after Jacobi, but Jacobi's work had soon a very
profound impact. For instance, one of the main contributions of Jacobi was the solution
of the {\bf inversion problem} explicitly for genus $g = 2$.

Concretely, Jacobi considered a polynomial $R(x)$ of degree $6$,
and then, given the two Abelian integrals
$$ u_1 (x_1, x_2) := \int_{x_0}^{x_1} \frac{dx} { \surd \overline{R}} +
\int_{x_0}^{x_2} \frac{dx} { \surd \overline{R}},$$
$$ u_2 (x_1, x_2) := \int_{x_0}^{x_1} \frac{x dx} { \surd \overline{R}} +
\int_{x_0}^{x_2} \frac{x dx} { \surd \overline{R}},$$
he found that the two symmetric functions
$$ s_1 := x_1 + x_2 ,  s_2 := x_1  x_2 ,$$
are $4$-tuple periodic functions of $u_1, u_2$.

Under the name {\bf Jacobi inversion problem} went the generalization
of this result for all genera $g$, and the solution to the Jacobi
inversion problem was one of the celebrated successes of Riemann.

Nowadays the result is formulated as follows: given a compact Riemann
surface $C$ of genus $g$, let $\omega_1, \dots \omega_g$ be a basis of
the space $H^0(\Omega^1_C)$ of holomorphic differentials on $C$
{\bf adapted} to a {\bf symplectic} basis 
$\alpha_1, \dots \alpha_g, \beta_1, \dots \beta_g$ for the Abelian
group of closed paths $H_1(C, \Z)$: this means that
$$ ( \int_{\alpha_i} \omega_j )= (\delta_{i,j}) = I_g ,
(\int_{\beta_i} \omega_j ) = (\tau_{i,j}),$$ where $I_g$ is the Identity
 $(g \times g)$ Matrix  and $\tau \in \mathcal H_g$,
 and that the intersection matrices satisfy $(\alpha_i , \alpha_j) = 0,
(\beta_i , \beta_j) = 0, (\alpha_i , \beta_j) = I_g.$
 
 Then the {\bf Abel-Jacobi} map $ C^g \rightarrow \C^g / (\Z^g + \tau
\Z^g) : = Jac (C)$, associating to the $g$-tuple $P_1, \dots P_g$ of
points of $C$ the sum of integrals (taken modulo $(\Z^g + \tau
\Z^g) $)
$$a_g(P_1, \dots P_g) := \int_{P_0}^{P_1} (\omega) + \dots
\int_{P_0}^{P_g} (\omega) ,$$ ($\omega$ being the vector with $i$-th component
$\omega_i$), 

is {\bf surjective} and yields a {\bf birational map} of the
{\bf symmetric product}
$$C^{(g)} := Sym^g (C) := C^g / \SSS _g,$$
$\SSS _g$ being the symmetric group of permutation of $g$
elements,  onto the {\bf Jacobian Variety of C} 

 $$ \ \ Jac (C) := \C^g /
(\Z^g + \tau \Z^g).$$

We see the occurrence of a matrix $\tau$ in the so called Siegel upper
half space $\mathcal H_g$ of symmetric matrices with positive definite
imaginary part: this positive definiteness ensures the convergence of
Riemann's theta series, and indeed Riemann used explicitly his theta
function to express explicitly the symmetric functions of the
coordinates of a $g$-tuple $P_1, \dots P_g$ of points of $C$ as
rational  functions of theta factors.

Today, we tend to forget about these explicit formulae, and we focus
our attention to the geometric description of the Abel-Jacobi maps for
any $n$-tuple of points of  $C$ in order to grasp the power of the
discoveries of Abel, Jacobi and Riemann.

The {\bf Abel-Jacobi} maps $ C^n \rightarrow \C^g / (\Z^g + \tau
\Z^g) : = Jac (C)$,  given by
$$a_n(P_1, \dots P_n) := \int_{P_0}^{P_1} (\omega) + \dots
\int_{P_0}^{P_n} (\omega) ,$$
and which naturally factor through the symmetric product $C^{(n)}$,
enjoy the following properties (cf. the  nice and concise lecture notes by D.
Mumford \cite{mum74} for more on this topic, and also \cite{corn} for a clear 
modern presentation) 

\begin{itemize}
\item
(1) {\bf Abel' s Theorem:} The fibres are projective spaces
corresponding to {\bf linearly equivalent} divisors, i.e. $ a_n(P_1,
\dots P_n) = a_n(Q_1, \dots Q_n)$ if and only if there is a rational
function $f$ on
$C$ with polar divisor  $ P_1 + \dots + P_n$  and divisor of zeros $
Q_1 +  \dots  + Q_n.$
\item
(2) For $n = g-1$  we have that the image of $a_{g-1}$ equals, up to
translation, the hypersurface $\Theta \subset Jac (C)$ whose
inverse image in $\C^g$ is given by the vanishing of the Riemann theta
function $\Theta= \{ z| \theta (z, \tau) = 0 \}.$
\item
(3) For $n = g$  we have that  $a_{g}$ is onto and birational.
\item
(4) For $n \geq 2g-1$  we have that  $a_{n}$ is a fibre bundle with fibres
projective spaces $\PP^{n-g}$.
\item
(5) For n=1 we have an embedding of the curve $C$ inside $Jac(C)$,
and in such a way that $C$ meets a general translate of the theta divisor
$\Theta$ in exactly $g$ points. And in this way one gets an explicit geometrical
description of the inverse to the Abel Jacobi map $a_g$. Indeed, to a point
$ y \in Jac (C)$ we associate the $g$-tuple of points of $C$ given by the intersection
$C \cap (\Theta + y)$ (here we think of $C$  $\subset Jac (C)$ under the embedding
$a_1$). 

\end{itemize}

Many of the above properties are very special and lead often to a characterization
of curves among algebraic varieties.

The above formulations are moreover the fruit of a very long process of maturation
whose evolution  is not easy to trace.
For instance, when did the concept of the {\bf Jacobian variety} of a curve
( Riemann surface) $C$ make its first appearance?

This notion certainly appears in the title of the papers by R. Torelli in 1913
(\cite{tor}) but apparently
\footnote{ We heard this claim from S.J. Patterson in G\"ottingen, soon after he
had written the article
\cite{pat99}.} the name was first used by F. Klein and became very soon extremely
popular. Observe however that in the classical treatise by Appell and Goursat
(\cite{appgour}), dedicated to analytical functions on Riemann surfaces, and
appeared first in 1895, although Jacobi's inversion theorem is amply discussed,
 no Jacobian or whatsoever variety is mentioned.

The French school of Humbert, Picard, Appell and Poincar\'e was very interested
about the study of the so called "hyperelliptic varieties", generalization of the
elliptic curves in the sense that they were defined as algebraic varieties
$X$ of dimension $n$ admitting a surjective entire 
holomorphic map \footnote{ The
letter $u$ clearly stands for "uniformization map".}
$ u: \C^n
\ra X$ .

Among those are the so called {\bf Abelian Varieties} \footnote{Which however,
 at the time of Torelli's paper, 1913, and also afterwards, were called Picard
Varieties . With the Prize winning Memoir by Lefschetz, \cite{lef}, the
terminology Abelian varieties became the only one in use.}, which are the
projective varieties which have a structure of an algebraic group.

In particular, Picard proved a very nice result in dimension $d=
2$, which was observed by Ciliberto (cf. his article in the present volume, also for
related historical references) to hold quite generally.
We want to give here a simple proof of this result, which we found during the
Conference, and which makes clear one basic aspect in which the higher dimensional
geometry has a different flavour than the theory of curves \footnote{ One could argue
whether the beauty of 1-dimensional geometry bears similarities to the surprising
isomorphisms of classical groups of small order. These also, by the way, are related
to Abel's heritage.
 For instance, the isomorphism $\SSS_4 \cong A(2, (\Z/2 ))$, whence $\SSS_3 \cong A(2,
(\Z/2 ))/(\Z/2 )^2  \cong PGL (2,(\Z/2 ))
\cong Aut (\HH)$. Here, $\HH$ is the group of order $8$ of unit integral quaternions,
and the last isomorphism is related to some later development in the theory of
algebraic curves, namely to Recilla's tetragonal construction, see \cite{rec},
\cite{don}}: namely there is no ramification in passing from Cartesian to symmetric
products.

\begin{teo}
Let $X$ be an  algebraic variety  of dimension $ d \geq 2$ and assume that there is
a natural number $n$ such that the $n$-th symmetric product $X^{(n)}$ is birational
to an Abelian variety $A$. Then $\ n=1$ (whence, $X$ is birational to an Abelian
variety).
\end{teo}

Moreover, this result ilustrates another main difference between the geometry
of curves and the one of higher dimensional varieties $X$: the latter has a quite
different flavour, because only seldom one can resort to the help of subsidiary
Abelian varieties for the investigation of a higher dimensional variety $X$.

\Proof
W.l.o.g. we may assume that $X$ be smooth.
Let us consider the projection $\pi$ of the Cartesian product onto the symmetric
product, and observe that $\pi : X^n \ra X^{(n)}$ is unramified in codimension $1$. 

We get a rational map $f :  X^n \ra  A$ by composing $\pi$ with the given birational
isomorphism, and then we observe that every rational map to an Abelian variety is a
morphism. 

It is moreover clear from the construction that $f$ is not branched in codimension
$1$: in particular,  it 
follows that
$X^n$ is  birational to an unramified covering of $A$, whence $X^n$ is  birational
to an Abelian variety.

Let us introduce now the following notation: for a smooth projective variety $Y$
we consider the algebra of global holomorphic forms  $H^0(\Omega_Y^*): =
\oplus_{i=0,
\dots dim (Y)} H^0(\Omega_Y^i) $ ({\bf holomorphic
algebra}, for short).

This graded algebra is a birational invariant, and for an Abelian variety $A$ it is the
free exterior algebra over $ H^0(\Omega_A^1) $. 

Now, the holomorphic algebra $H^0(\Omega_{X^n}^*)$ of a Cartesian product $X^n$ is
the tensor product of $n$ copies of the holomorphic algebra $H^0(\Omega_{X}^*)$ of $X$.

Denote $H^0(\Omega_{X}^*)$ by $B$: then we reached the conclusion that $B^{\otimes n}$
is a free exterior algebra over its part of degree $1$, $B_1^1 \oplus B_2^1  \dots
\oplus B_n^1$.

It follows that also $B$ is a free exterior algebra (i.e.,
$X$ enjoys the property that its holomorphic algebra
$H^0(\Omega_{X}^*)$ is a free exterior algebra over $H^0(\Omega_{X}^1)$).

Moreover, the holomorphic algebra of the symmetric product is the 
invariant part of this
tensor product (for the natural action of the symmetric group in $n$-letters
$\SSS_n$), and by our assumption  $C:= (B^{\otimes n})^{\SSS_n} $ is also a free
exterior algebra. However, $C^1 = B^1$, and $r:= dim (B^1)$ is also the highest
degree $i$ such that $ B^i \neq 0$. But then $C^{nr} \neq 0$, contradicting the
property that $C$ is a free exterior algebra with $ dim \ (C^1) = r$, if $n \neq 1$. 

Thus $n=1$, and $X$ is  birational to an Abelian variety.

\qed

\begin{oss}
In the case of dimension $d=1$, the same algebraic arguments easily yield that the
n-th symmetric product of a curve $C$ is not birational to an Abelian variety if $n
\neq g$, $ g:= dim \ (H^0(\Omega^1_C))$.
\end{oss}

Once more, the algebra of differential forms, as in Abel's work, has played the
pivotal role.

The importance of this algebra was observed also by Mumford (\cite{mum68}) who
used it to show that on an algebraic surface $X$ with $H^0(\Omega^2_X) \neq
0$, the group of $0$-cycles ( Sums $\Sigma _im_i P_i$ of points $P_i \in X$ with
integer multiplicities $m_i \in \Z$) modulo rational equivalence is not
finite dimensional, contrary to the hope of Severi, ( D. Mumford  sarcastically
wrote: "One must admit that in this case the {\em technique} of 
the italians was
superior to their vaunted intuition" \footnote{ However, as well known, there are
Italians with techniques and  ideas, and others who are not perfect. In particular,
while it is not difficuult to find errors or wrong assertions in Enriques and
Severi, it is rather hard to do this with Castelnuovo. }) who unfortunately was
basing his proposed theory on a wrong article  (\cite{sev32}, where not by chance
the error was an error of ramification).

It must be again said that the 
italian school, and especially Castelnuovo, gave a
remarkable impetus to the geometrization of the theory of Abelian varieties.

This approach, especially through the work of Severi, influenced Andre' Weil who
understood the fundamental role of Abelian varieties for many questions of 
algebraic number theory.  Weil used these ideas to construct (\cite{weilva}) the
Jacobian variety of a  curve as a quotient of the symmetric product
$C^{(g)}$, and then, for a $d$-dimensional variety with $d \geq 2$, the Albanese
variety $Alb(X)$ as a quotient (in the  category of Abelian varieties) of the 
Jacobian $J(C)$ of a sufficiently general linear section $ C = X \cap H_1 \cap H_2
\dots \cap H_{d-1}$.

It must be however said that also the later geometric constructions 
were deeply
influenced by the bilinear relations which Riemann, through a convenient dissection
of his Riemann surface $C$, showed to hold for the periods of the Abelian integrals
of  the first kind of $C$.

Nowadays, the usual formulation is (according to Auslander and Tolimieri,
\cite{a-t} pages 267 and 274, the first formulation is essentually due to Gaetano
Scorza in
\cite{scor16}, while the second is essentially due to Hermann Weyl in
\cite{Weyl34}, \cite{Weyl36}, with refinements from A. Weil's book \cite{weilVK})

\begin{df}
Let $\Gamma$ be a discrete subgroup of a complex vector space $V$, such that the
quotient $V/\Gamma$ is compact (equivalently, $\Gamma \otimes \R \cong V$): then
we say that the complex torus  $V/\Gamma$ satisfies the  two Riemann bilinear
relations if

\begin{itemize}
\item
I) There exists an alternating form $ A : \Gamma \times \Gamma \ra \Z$
such that $A$ is the imaginary part of an Hermitian form $H$ on $V$
\item
II) $H$ is positive definite.

\end{itemize}
\end{df}
\begin{oss}
Or, alternatively, a complex structure on $\Gamma \otimes \C $, i.e., a
decomposition $\Gamma \otimes \C  = V \oplus \bar{V}$ and an element  
$ A \in \Lambda^2 (\Gamma)^{\vee}$ yield a polarized Abelian variety 
if the component of $A$ in $\Lambda^2 (V)^{\vee} \subset \Lambda^2
(\Gamma \otimes \C )^{\vee}$ is  zero and then its component 
in $(V)^{\vee} \otimes (\bar{V})^{\vee}$ is a positive definite Hermitian form.
\footnote{ The second characterization is very useful for the study of
fibre bundles of Abelian varieties, as we had opportunity to experience
ourselves, cf.
\cite{cat02-2} }
\end{oss}

The basic theorem characterizing complex Abelian varieties is however due to
Henri Poincar\'e (\cite{poi84}, \cite{poi02}) who proved the linearization of the
system of exponents, 
i.e., the more difficult necessary condition in the theorem, by an averaging
procedure (integrating the ambient Hermitian metric of $ X \ra
\PP^N$ with respect to the translation invariant measure of  $X = V/\Gamma$, he
obtained a translation invariant Hermitian metric).

\begin{teo}
A complex torus $X = V/\Gamma$ is an algebraic variety if and only if the two
Riemann bilinear relations hold true for $V/\Gamma$.

Both conditions are equivalent to the existence of a meromorphic function $f$ on
the complex 
vector space $V$ whose group of periods is exactly $\Gamma$
(i.e., $\Gamma = \{ v \in V| f(z + v) \equiv f (z) \} ). $
\end{teo}

Poincar\'e had an extensive letter exchange with Klein (cf. Klein's Collected
Works, where pages 587 to 621 of Vol. III are devoted to the "Briefwechsel" between
the two, concerning the problem of uniformization, and their early attempts, which
were based on a 'principle of continuity' which was not so easy to justify \footnote{
In his unpublished Fermi Lectures held in Pisa in 1976, D. Mumford explained how this
approach was working, using clarifications due to Chabauty, \cite{chab} }), especially
related to the study of discontinuous groups, acting not only on
$\C^n$ as in the case of tori, but also on the hyperbolic upperhalf 
plane $\HH$. The
main result,  whose complete proof was obtained in 1907 by a student of Klein, Koebe,
and by Poincar\'e independently, was the famous uniformization theorem that again we
state in 
its modern formulation for the  sake of brevity. 

\begin{teo}
If a Riemann surface is not the projective line $\PP^1_{\C}$,  the complex
plane $\C$, nor $\C^*$ or an elliptic curve, then its universal covering is the
(Poincar\'e) upper half plane $\HH$.
\end{teo}

The reason why the upper half plane "belongs" to Poincar\'e is that Klein
prefered to work with the biholomorphically equivalent model given by the unit
disk $\D := \{ z \in \C | |z| < 1 \}.$ In this way Klein was capable of making
us the gift of beautiful symmetries given by tesselations of the disk by
fundamental domains for the action of very explicit Fuchsian groups (discrete
subgroups 
$\Gamma$ of $PSU(1,1,\C)$) with compact quotient $\D / \Gamma$.

To summarize the highlights of the turn of the century, 
when geometry was a very
central topic, one should say that  several new
geometries came to birth at that time: but  the new developments were
based on new powerful analytic tools, which were the bricks of the new building.

However, although the birth of differential geometry lead to new geometrical
theories based on infinite processes 
where metric notions played a fundamental
role, algebraic geometry went on with alternating balance between 
geometrical versus algebraic methods.

\section{Algebraization of the geometry}

At this moment, a  witty reader, tired of the distinction "algebra"-"non
algebra", might also 
remind us that  the popular expression "This is algebra
for me" simply means: 'I do not understand a single word of this'.

There is a serious point to it: the concept algebra is slightly ambiguous, and a
very short formula could  be not very inspiring without a thoroughful explanation
of its meaning(s), and of all the possible consequences and applications.

One of the best ways to understand a formula is for instance to relate it to a
picture, to see it thus related to a geometrical or dynamical process.

Needless to say, the best example of such an association is the Weierstrass 
equation of a plane cubic curve 
$$ C^1_3 = \{ (x,y,z) \in \PP^2_{\C} | y^2 z = 4 x^3 - g_2 x z^2 - g_3 z^3 \}. $$
To this equation we immediately associate the picture yielding the group law of
$C=C^1_3$, i.e., a line $L$ intersecting $C$ in the three points $P, Q$ and $T=
- (P+Q)$. 

Where have we seen this picture first?\footnote{ Not always mathematical concepts have
a birthday. But sometimes it happens , as for another Legacy by Abel: the concept of an
abstract group was born in 1878, with a "Desiderata and suggestions" by A. Cayley
\cite{cay1}. In this note appears first the multipication table of a group.
Immediately afterwards, however, in \cite{cay2}, Cayley realized 
that it is much better to
work with what is now called the "Cayley graph" of a group 
endowed with a set of generators. }
Well, in my case, I (almost) saw it first in the book by Walker on Algebraic curves,
exactly in the last paragraph, in the section 9.1 entitled "Additions of points on a
cubic". The book was written in 1949, and if we look at books on algebraic curves
written long before, the group law is not mentioned there. For instance,  Coolidge's
book "A treatise on algebraic plane curves" has a paragraph entitled "elliptic curves"
, pages 302-304, and the main theorems are first that an elliptic curve is birational
to a plane cubic, and then the {\bf Cross-Ratio Theorem} asserting that if $P\in C^1_3
\subset
\PP^2$ is any point, through $P$ pass exactly $4$ tangent lines, and their cross ratio
is independent of the choice of $P \in C$.
The Weierstrass equation, and the explanation that the cubic curve is
uniformized through the triple $(1,\PPP, \PPP')$, where $\PPP$ is the
Weierstrass function, comes later, as due after the Riemann Roch theorem
\footnote { 
"With proof communicated to the author by prof. Osgood verbally in
Nov. 1927".}, on pages 363-367 in the paragraph "Curves of genus 1".

Going to  important textbooks of the Italian tradition, like Enriques and
Chisini's 4 Volumes on the "Lezioni sulla teoria geometrica delle equazioni e
delle funzioni algebriche"\footnote{ Edited by Zanichelli in the respective
years 1915, 1918, 1924, 1934.} we see that Volume IV contains the Book 6,
devoted to "Funzioni ellittiche ed Abeliane". Here Abel's theorem is fully
explained, and on page 77 we see Abel's theorem for elliptic curves, on page 81
the addition theorem for the $\PPP$ function of Weierstrass: the geometry of the
situation is fully explained, i.e., that three collinear points sum to zero in
the group law given by the sum of Abelian integrals of the first kind
$$ u_1 + u_2 + u_3 = \int_{x_0}^{x_1} \frac{dx}{y} + \int_{x_0}^{x_2}
\frac{dx}{y} +\int_{x_0}^{x_3} \frac{dx}{y} = 0.$$
It is also observed that the inverse of the point $P = (1,x,y)$ is the point
$P = (1,x,-y)$, as a consequence of the fact that $\PPP$ is an even function.

This Book 6 is clearly influenced by Bianchi's Lecture 
Notes " Lezioni sulla
teoria delle Funzioni di variabile complessa"\footnote{ Spoerri, Pisa, 1916}
whose Part 2 is entirely devoted to "Teoria delle funzioni ellittiche", and in
the pages 315-322 the addition theorem for $\PPP$ is clearly explained, moreover
"Alcune applicazioni geometriche" are given in the later pages 415-418.

So, the picture is there, is however missing the wording: a plane cubic is an Abelian
group through the sum obtained via linear equivalence of divisors, namely, the sum of
three points $P,Q,T$ is zero if and only if the divisor $P+Q+T$ is linearly
equivalent to a fixed divisor $D$ of degree 3 (in the Weierstrass model, $D$ is
the divisor $3O$, $O$ being the flex point at infinity).

As it was explained to me by Norbert Schappacher (cf. \cite{schap90}), the works of
Mordell and Weil in the 1920's are responsible of this new wording and perspective.
In fact, these authors considered a cubic curve $C$ whose equations has
coefficients in a field $K$, and noticed that the set $C(K)$ of
$K$-rational points, i.e., the points whose coordinates are in $K$, do indeed
form a subgroup. For this they did not need that one flex point should be
$K$-rational, since essentially, once we have a $K$-rational point $O$, we can
reembed the elliptic curve $C$ by the linear system $|3O|$, and then obtain a
new cubic $C'$ whose $K$ rational points are exactly those of $C$.

Through these works started an exciting new development, namely the
geometrization of arithmetic, which was one of the central developments in the
20-th century mathematics.

For instance, the theory of elliptic curves over fields of finite characteristic
was (cf\cite{weil29}) built by Weil and then Tate who (cf. the quite late
appearing in print \cite{tate}), starting from the Weierstrass equation, slightly
modified into
$$y^2 z = x^3 - p_2 x z^2 - p_3 z^3$$
started to construct analogues of the theory of periods.\footnote{De hoc satis,
because the talks by Faltings and Wiles were exactly dealing with these aspects
of Abel's legacy.}

Going back to Bianchi, it is Klein's and Bianchi's merit to have popularized the
geometric picture of elliptic curves, and actually Bianchi went all the way
through in some of his papers to describe the beautiful geometry related to the
embeddings of elliptic curves as non degenerate curves of degree $n$ in
$\PP^{n-1}$ ($\forall n \geq 3$).

It took however quite long till a purely algebraic interpretation of
Weierstrass' equation made its way through.

Nowadays we would associate to an elliptic curve $C$ and to a point $O \in C$
the $\N$-graded ring 
$$ \bigoplus_{m=0}^{\infty} H^0 (C, \hol_C(mO)): = \RR (C, \LL)  $$
where $C$ is defined over a (non algebraically closed) field $K$ of characteristic
$\neq 2$, 
$O$ is a $K$-rational point and $\LL :=\hol_C(O)$. 

As a consequence of the Riemann Roch theorem we obtain the following

\begin{teo}
$$\RR (C, \LL)  \cong  K [u,\xi, \eta ] / ( \eta^2 - \xi^3 + p_2 \xi u^4 + p_3
u^6 ),$$
where $deg(u)=1$, $deg(\xi)=2$,$deg(\eta)=3$, $div(u) = O$.
\end{teo}

To go back to the 
original Weierstrass equations it suffices to observe that
$$ \PPP = \xi / u^2, \PPP' = \eta/u^3 ,$$ and that  
$ x := \xi  u, y := \eta , z := u^3 $ are a basis of the vector space $H^0 (C,
\hol_C(3O))$.

One sees also clearly how the Laurent expansion at $O$ of $\PPP$ is
determined by and determines $p_2, p_3$.

Surprisingly, the following general problem is still almost completely open,
in spite of a lot of research in this or similar directions

\begin{problem}
Describe the graded ring $\RR (A, \LL)$ for $\LL$ an ample divisor on an Abelian
variety, for instance in the case where $\LL$ yields a principal polarization.
\end{problem}

Before explaining the status of the question, I would first like to explain its
importance.

Take for instance the case of an elliptic curve $C$ whose ring is completely
described (the ring does not depend upon the choice of $O$ because we have a
transitive group of automorphisms provided by translations for the group law of
$C$).

We want for instance to describe the geometry  of the embedding of $C$ as a curve
of degree $4$ in $\PP^3$. We observe that, at least in the case where $K$ is
algebraically closed, any such embedding is given by the linear system $|4 O|$,
for a suitable choice of $O$.

The coordinates of the map are given by a basis of the vector space $H^0(C,
\hol_C(4O))$, i.e., by 4 independent homogeneous elements of degree $4$ in our
graded ring $\RR (C, \LL)$. These are easily found to be equal to
$s_0 := u^4, s_1 := u^2 \xi, s_2 :=\xi^2, s_3 :=u \eta$. Then we obviously have
the two equations 
$$ s_0 s_2 = s_1 ^2 , \ \ s_3^2 =  s_1 s_2 + p_2 s_0 s_1 + p_3
s_0^2 ,$$   holding for the image of $C$ (the second is obtained by the
"Weierstrass" equation once we multiply by $u^2$).

These are all the equations, essentially by Bezout's theorem, since
$C$ maps to a curve of degree $4$ and $ 4 = 2 \times 2$.

From an algebraic point of view, what we have shown is the process of
determining a subring of a given ring, and the nowadays computer algebra
programs like "Macaulay" have standard commands for this operation
(even if sometimes the computational complexity of the process may become too
large if one does not use 
appropriate tricks).

Classically, a lot of attention was devoted to the geometric study of the maps
associated to the linear systems $| m \Theta|$, where $\Theta$ is the divisor
yielding a principal polarization of the given Abelian variety. 

For instance, the $4$-ic Kummer surface is the image in $\PP^3$ of a principally
polarized Abelian surface under the linear system 
$| 2 \Theta|$, yielding a $ 2 : 1$ morphism which identifies a point $v$ to
$-v$, and blows up the $16$ $2$-torsion points to the $16$ nodal singularities
of the image surface.

There is a wealth of similar results, which can be found for instance in the 
books of Krazer \footnote{"Theorie der Thetafunktionen, Teubner, 1894.},
Krazer-Wirtinger, and Coble (\cite{kraz}, \cite{k-w}, \cite{coble}), written in
the period 1890-1926. 

The first books are directly influenced by the Riemann quadratic relations,
i.e., linear relations between degree two monomials in theta functions with
characteristics (and their coefficients being also products of "Thetanullwerte",
i.e., values in $0$ of such thetas with characteristics), and show an attempt to
use 
geometrical methods starting from analytic identities. 
Coble's book is
entitled "Algebraic geometry and theta functions", and is already influenced by
the breakthrough made by Lefschetz in his important Memoir (\cite{lef}).

Lefschetz used systematically the group law to show that if $s_1(v), \dots
s_r(v) \in H^0(m \Theta)$ and we choose points $a_1, \dots a_r $ such that 
$ \Sigma_i a_i = 0$, then the product $s_1(v+a_1) \dots
s_r(v + a_r) \in H^0(r m \Theta)$. Then he chooses sufficiently many and
sufficiently general points $a_1, \dots a_r$  so that these sections separate points
and tangent vectors. Thus Lefschetz proves in particular

\begin{teo}
$| m \Theta|$ yields a morphism for $ m\geq 2$ and an embedding for $m \geq 3$.
\end{teo}

The direction started by Lefschetz was continued by many authors, notably Igusa,
Mumford, Koizumi, Kempf, who proved several results concerning the equations of
the image of an Abelian variety (e.g., that the image of $| m \Theta|$  is an
intersection of quadrics for $m \geq 5$).

I will later return  on some new ideas related to these developments,  I
would like now to focus on the status of the problem I mentioned.

The case where $A$ has dimension $g=1$ being essentially solved by Weierstrass,
the next question is whether the answer is known for $g=2$. This is the case,
since the description of the graded ring was obtained by A.
Canonaco in 2001 (\cite{can02}); an abridged version of his result is as follows

\begin{teo}
Let $A$ be an Abelian surface and $\Theta$ 
be an effective divisor yielding a
principal polarization: then the graded ring $\RR (A, \hol_A(\Theta))$ has a
presentation with $11$ generators, in degrees $ (1,2^3,3^5,4^2)$, and $37$
relations, in degrees $ (4,5^6,6^{17},7^{10}, 8^3)$.
\end{teo}

\begin{oss}
1) Canonaco gives indeed explicitly the $37$ relations, whose shape however is
not always the same. One obtains in this way an interesting stratification of the
Moduli space of p.p. Abelian surfaces. Does this stratification have a simple
geometrical meaning in terms of invariant theory?
\end{oss}

The proof uses at a certain point some computer algebra aid, since the equations
are rather complicated.

Nevertheless, we would like to sketch the simple geometric ideas underlying the
algebraic calculations, since, as one can easily surmise, they are
related to the aforementioned Riemann's developments of Abel's investigations.

\Proof
The key point is thus that, in the case where $\Theta$ is an irreducible
divisor,  $\Theta$ is isomorphic to a smooth curve $C$ of genus $2$ (the other
case where $\Theta$ is reducible is easier, since then $A$ is a product of
elliptic curves, and $\Theta$ is the union of a vertical and of a horizontal
curve).

One uses first of all the exact sequence (for $n\geq 2$, 
since for $n=1$ the right
arrow is no longer surjective)
$$0 \ra H^0(A,\hol_A( (n-1)\Theta)) \ra  H^0(A,\hol_A( n\Theta)) \ra
H^0(\Theta,\hol_{\Theta}( n\Theta))  \ra 0$$ and of the isomorphism
$$ H^0(\Theta,\hol_{\Theta}( n\Theta)) \cong H^0(C,\hol_C( n K_C)).$$
One relates thus our graded ring to the canonical ring of the curve $C$,
which is well known, the canonical map of $C$ yielding  a double covering of
$\PP^1$ branched on $6$ points.

In more algebraic terms, there is a homogeneous polynomial $R(y_0,y_1)$ of
degree $6$ such that
$$ \RR (C,\hol_C(  K_C) ) \cong K [y_0,y_1,z] / (z^2 - R(y_0,y_1)).$$

One can summarize the situation by observing that, if $ \RR_A := \RR (A,
\hol_A(\Theta))$, $\RR_C := \RR (C,\hol_C(  K_C) )$, then $\RR_A$ surjects onto
the subring $\RR'$ of $\RR_C$ defined by $\RR' : = \bigoplus_{n\geq 2}
H^0(C,\hol_C( n K_C))$.

To lift the ring structure of $\RR'$, which is not difficult to obtain, to the
ring structure of $\RR_A$ we use again Abel's theorem, i.e. the sequence of maps
$$ C^2 \ra C^{(2)} \ra A  $$
where the last is a birational morphism contracting to a point the divisor $E$
consisting of the set of pairs $\{ (P, i(P))| P \in C \}$ in the canonical
system of $C$ (here $i$ is the canonical involution of the curve $C$).

Letting $D_i$ be the pull back of a fixed Weierstrass point of $C$ under the
$i$-th projection of $C \times C$ onto $C$, we obtain that
$$  H^0(A,\hol_A( n\Theta)) \cong  H^0(C^2,\hol_{C^2}( n (D_1 + D_2 +E)))^+, $$ 
where the superscript $^+$ denotes the $+1$-eigenspace for the involution of
$C^2$ given by the permutation exchange of coordinates. It also helps to consider
that the
$\Theta$ divisor of $A$ is the image of a vertical divisor $ \{P\} \times C$.
Instead, the diagonal $\Delta_C$ of $C^2$ enters also in the picture because the
pull back (under $h= \phi_K \times \phi_K : C \times C \ra \PP^1 \times \PP^1$,
$\phi_K$ being the canonical map of $C$) of the diagonal of $\PP^1 \times \PP^1$
is exactly the divisor $\Delta_C + E$, whence $\Delta_C + E$ is linearly
equivalent to $2 D_1 + 2 D_2$.

We omit the more delicate parts of the proof, which are however based on the
above linear equivalences and on the  action of the dihedral group $ D_{4\times
2}$ on
$C \times C$ (this is a lift, via the $(\Z/2)^2$ Galois cover  $h:  C \times C
\ra
\PP^1 \times \PP^1$ of the permutation exchange of coordinates on  $ \PP^1
\times \PP^1$).

\qed

\begin{oss}
1) This approach should work in principle for the more general case of  the
Jacobian variety of a hyperelliptic curve. In fact,  for each Jacobian variety
we have the sequence of maps
$$ C^g \ra C^{(g)} \ra A  $$
and the $\Theta$ divisor is the image of a big vertical divisor $ \{P\} \times
C^{g-1}$. Again we have a $(\Z/2)^g$ Galois cover $h : C^g \ra (\PP^1)^g$ and a
semidirect product of the Galois group with the symmetric group in $g$ letters
$\SSS_g$ 
(the group of the $g$-dimensional cube).

2) Another question is whether there does exist a more elegant, or just shorter
presentation for the ring.

\end{oss}

{\bf  THE ALGEBRA OF THETA FUNCTIONS AND THE ALGEBRA OF REPRESENTATION THEORY} 

A more conceptual understanding of the several identities of general theta
functions came through the work of Mumford (\cite{mum66-7} , cf. also \cite{igu}
and \cite{weil64}).

In  Mumford's articles and in Igusa's treatise one finds  a clear path set by choosing
representation theory as a guide line, especially as developed by Weyl, Heisenberg and
von Neumann.

The basic idea is shortly said: let $G$ be a compact topological group, endowed
henceforth with the (translation invariant) Haar measure $d \mu_G$. 

Consider then the vector space $V= L^2(G, \C)$: then we have an action $\tau$ 
of the group $G$ on $V$ by translations, $\tau_{\gamma} (f )(g) := f (g \gamma
^{-1})$.

Defining the group of characters $G^*$ as $\ G^* := Hom (G, \C^*)$ we have an
action of $G^*$ on $V$ given by multiplication $ \chi f (g) := \chi(g) f(g)$.

The two actions fail to commute, but by very little, since 
$$ \chi [\tau_{\gamma} (f )] (g) = \chi(g)f (g \gamma ^{-1}) $$ $$
 \tau_{\gamma} ( \chi f ) (g) = \chi(g \gamma ^{-1})f (g \gamma ^{-1}) =
 \chi( \gamma) ^{-1} \chi (g)f (g \gamma ^{-1}) $$
 thus commutation fails just up to multiplication with the constant function
$\chi( \gamma) ^{-1}$.

Together, the action of $G$ and of $G^*$ generate a subgroup of the Heisenberg group,
a central extension 
$$ 1 \ra \C^* \ra Heis (G) \ra G \times G^* \ra 1.$$
It turns out that the algebra of theta functions is deeply related to the
representation theory of the Heisenberg group of the Abelian variety $A$
(we see $A$ as the given group
$G$). 

But, as Mumford pointed out, we have a more precise relation which takes into
account a given line bundle $\LL$.

In the case of an Abelian variety $A$,  the group of characters is endowed with
a complex structure viewing it as the {\bf Picard variety} $Pic^0(A)$, the connected
component of $0$ in$H^1(A,\hol_A^*)$. $Pic^0(A)$ is also called the dual Abelian
variety, and a non degenerate line bundle $\LL$ is one for which the homomorphism
$\phi_{\LL}  : A \ra Pic^0(A)$, defined by $\phi_{\LL} (x) = T_x^*(\LL) \otimes
\LL^{-1}$ 
($T_x$ denoting translation by $x$), is surjective (hence with finite
kernel $K(\LL)$).

Mumford introduces the finite Heisenberg group associated to $\LL$ 
via the 
so called Thetagroup of $\LL$, defined as ${\rm Theta}(\LL) := \{ (x, \psi)|
\psi:
\LL
\cong T_x^* (\LL) \}$.

${\rm Theta}(\LL)$ is a central extension of $K(\LL)$ by $\C^*$, but since
$K(\LL)$ is finite, if $n$ is the exponent of $K(\LL)$, the central extension is
induced (through extension of scalars) by another  central extension
$$ 1 \ra \mu_n \ra \Theta(\LL) \ra K(\LL) \ra 1,$$
where $\mu_n$ is the group of $n$-th roots of unity.

Moreover, the alternating form $\alpha : \Gamma \times \Gamma \ra \Z$ given by
the Chern class of $\LL$ gives a non degenerate symplectic form 
on $K(\LL)$ with
values in $\mu_n$, thus allowing to easily obtain from $\Theta(\LL)$  the Heisenberg
group of a finite group $G$.

The geometry of the situation is that the group $K(\LL)$ acts on the projective
space associated to the vector space $H^0(A, \LL)$: but if we want a linear
representation on the vector space
$H^0(A, \LL)$ we must see this vector space as a representation of the
finite Heisenberg
 group $\Theta(\LL)$ (thus we have a link with Schur's theory of
multipliers of a projective representation).

In the case of $\LL = \hol_A (n \Theta)$, $K(\LL)$ consists of the
subgroup $A_n$ of n-torsion points, and another central idea, when we have to
deal with a field of positive characteristic, is to replace the vector space
$\Gamma \otimes \R$ with the inverse limit of the subgroups $A_n$.

The story is too long and too recent to be further told here: Mumford used this
idea in order to study the Moduli space of Abelian varieties over fields of
positive characteristics, and in turn this was used to take the reduction
modulo primes of Abelian varieties defined over number fields.

These results were crucial for arithmetic applications, especially Faltings'
solution (\cite{falt83}) of the

{\bf Mordell conjecture. Let a curve $C$  of genus $g\geq 2$ be defined over a
number field $K$: then the set $C(K)$ of its rational points is finite.} 

It must be furthermore said that the algebraic calculations allowed by the study of the
characters of representations of the finite Heisenberg groups has lead
also to a better concrete understanding of equations and geometry of Abelian varieties.

Surprisingly enough, even in the 
case of elliptic curves this has led, together with
Atiyah's study of vector bundles on elliptic curves (cf.\cite{ati57}) to a deeper
understanding of the geometry of symmetric products of elliptic curves and their 
maps to projective spaces (cf. \cite{ca-ci}).

The recent literature is so vast that we have chosen to mention just a single but
quite beautiful example, due to Manolache and Schreyer (cf. \cite{man-schr}).

The authors give several equivalent descriptions of the moduli space $X(1,7)$ of
Abelian surfaces $S$ with a polarization $L$ whose elementary divisors are
$(1,7)$.

Their main result is that this moduli space is birational to the Fano 3-fold
$V_{22}$ of polar hexagons to the Klein plane quartic curve $C$ ( of equation $ x^3 y +
y^3 z + z^3 x = 0 $) which is a compactification of the moduli space $X(7)$ of elliptic
curves $E$ with a level $7$ structure, i.e., elliptic curves given with an
additional isomorphism of the group $E_7$ of torsion points with $(\Z / 7)^2$.

The Klein quartic is rightly famous because it admits then as group of
automorphisms the group
$\PP SL (2,\Z / 7) = SL (2,\Z / 7) / \{ \pm I \}$, a group of cardinality $168 $, and
as it is well known this makes the Klein quartic the curve of genus $3$ with the
maximal number of automorphisms (cf. \cite{acc}).

Now, it is easy to suspect  some connection between the pairs $(S,L)$ and the
pairs $(E, \LL : = \hol_E(7 O))$ once we have learnt of the finite Heisenberg
group: in fact, the 
respective groups for $L$ and for $\LL$ are isomorphic, and the
respective complete linear systems $|L|$ and $|\LL|$ yield  embeddings for
$E$, respectively for the general $S$, into $\PP^6$.

That is, we view both $E$ and $S$ as Heisenberg invariant subvarieties of the
same $\PP^6$ with an action of $(\Z/7)^2$  provided by the projectivization of
the standard representation of the Heisenberg group on $\C^{\Z/7}$.

The geometry of the situation tells us that $E$ is an intersection of $28 - 14 =
14$ independent quadrics, while we expect $S$ to be contained in $28-  28= 0$
quadrics, so we seem to be stuck without a new idea.

The central idea of the authors is to think completely in algebraic terms,
looking at a self dual locally free Hilbert resolution of the ideal of $S$,
which has length $5$ instead of $4$ (because of $H^1(S, \hol_S) \neq 0$).

It turns out that the middle matrix, because of Heisenberg symmetry, boils down
to a $3 \times 2$ matrix of linear forms on a certain $\PP^3$. Then, the $3$
determinants of the $2 \times 2$-minors yield three quadric surfaces whose
intersection is a twisted cubic curve $\Gamma_S$ which is shown to completely
determine $S$. In this way one realizes the moduli space as a certain subvariety
 of the Grassmann variety of $3$-dimensional vector subspaces of the
vector space $H^0(\PP^3, \hol_{ \PP^3 } (2))$ of quadrics in the given $\PP^3$
(one takes the locus of subspaces where a certain antisymmetric bilinear map
 of vector bundles restricts to zero).

This subvariety is the Fano $3$-fold $V_{22}$ mentioned above, which is a rational
variety. 

At this point we don't want to deprive the reader of the pleasure of learning
the intricate details  from the original sources (\cite{man-schr}, \cite{schr}):
but we need at least to explain what is a polar hexagon of a plane quartic curve
$C$ with equation $ f(x) = 0$. 

Observe that the polynomial $f(x)$ depends upon exactly 
$15$ coefficients,
while, if we take $6$ linear forms $l_i(x)$, they depend on 
$18$ affine parameters,
and we expect therefore to have a $3$-dimensional variety $Hex(C)$ parametrizing the
$6$-tuples of such linear forms such that $f(x) = \Sigma_{i} l_i(x)^4.$

$Hex(C)$ is called the variety of polar hexagons, and it is indeed a $3$-fold
in the case of the Klein quartic. Finally, these constructions allow to find the Klein
quartic as the discriminant of the net of quadrics in $\PP^3$ associated
to the Hilbert resolution of $S$.

Also the identification of the Fano $3$-fold with $Hex(C)$ is based on the
study of higher syzygies, but beyond this many other beautiful classical results are
used, which are due to Klein, Scorza, and Mukai (cf. \cite{klein78}, \cite{scor99},
\cite{muk92}).

Especially nice is the old theorem of Scorza (\cite{scor99}), proved in 1899,
that the variety of plane quartic curves is birational to the variety of
pairs of a plane quartic curve $D$ given together with an even
theta-characteristic (this amounts to 
writing the equation of $D$ as the
discriminant of a net of quadrics in $\PP^3$). This theorem is a clear  example of the
geometrization of ideas coming from the theory of theta functions (which, as we
saw, are certain Fourier series, and therefore, seemingly, purely analytical
objects).

\bigskip

{\bf ABELIAN VARIETIES AND MULTILINEAR ALGEBRA}

The roots of these developments, which historically go under the name "The problem
of Riemann matrices", and occupied an important role for the birth of the theory
of rings, modules and algebras, are readily explained by the following basic

\begin{oss}
Given two tori $T = V/ \Gamma, T' = V'/ \Gamma'$, any holomorphic map 
$f : T \ra T'$ 
between them is induced by a complex linear map $ F: V \ra V'$ such that $ F
(\Gamma) \subset \Gamma'$.
\end{oss}

Whence, for a $g$-dimensional torus, the ring 
$$End(T) 
:= \{ f : T \ra T | f \ {\rm is \ holomorphic }, f(0) = 0\}  $$ is the subring
 of the ring 
 of matrices $ Mat (2g, 2g, \Z) $ given by
 $$End(T) :\{ B \in End ( \Gamma) | B \in   (V^{\vee} \otimes V ) \oplus
(\bar{V}^{\vee}
\otimes \bar{V})  \subset 
 (\Gamma \otimes \C)^{\vee} \otimes (\Gamma \otimes \C)\}  $$
 since then the restriction of $B $ to $V \subset (\Gamma \otimes \C)$ is complex
linear.

In general, the study of endomorphism rings of complex tori is not completely
achieved.

The main tool which makes the case of Abelian Varieties easier is the
famous {\bf Poincar\'e's complete reducibility theorem} (cf. \cite{poi84})
\begin{teo}
Let $A'$ be a subabelian variety of an Abelian variety $A$: then there exists
another Abelian variety $A"$ and an {\bf isogeny} (a surjective homomorphism with
finite kernel) $ A' \times A" \ra A$.
\end{teo}
\Proof
The datum of $A'$ amounts to the datum of a sublattice $\Gamma' \subset \Gamma $
which is saturated ( $ \Gamma /\Gamma'$ is torsion free) and complex (i.e., there
is a complex subspace $W \subset V$ with 
$ W \oplus  \overline{W 
} = \Gamma' \otimes \R \subset  \Gamma \otimes \R = V \oplus  \overline{V } $ ).

Now, given the alternating form $A$, its orthogonal in $\Gamma$ yields a
sub-lattice $\Gamma "$ spanning the complex subspace $U$ orthogonal to $W$ for the
Hermitian bilinear product associated to $H$: since $H$ is positive definite,
we obtain an orthogonal  direct sum $ V = W \oplus U$, and we define
$ A" : = U /\Gamma"$.

\qed

We sketched the above proof just with the purpose of showing how   the
language of modern multilinear algebra is indeed very appropriate for these types
of questions.

The meaning of the reducibility theorem is that, while for general tori a subtorus
$ T' \subset T$ only yields a quotient torus $ T / T'$, here we get a direct sum
if we consider an equivalence relation which identifies two isogenous Abelian
vareties.

Algebraically, the winning trick was thus to classify first Endomorphism Rings
tensored with the rational integers, because

\begin{oss}
If $T = V / \Gamma, T' = V / \Gamma'$ are isogenous tori, then $End(T) \otimes \Q
\equiv End(T') \otimes \Q$.
\end{oss}

And then the study is restricted to the one of {\bf Simple Abelian Varieties},
i.e., of the ones which do not admit any Abelian subvariety whatsoever (naturally,
this concept was very much inspired by the analogous concept of curves which do
not admit a surjective and not bijective mapping onto a curve of positive genus).

The classification of Endomorphism rings of Abelian varieties was achieved through
a long series of works by Scorza, Rosati, Lefschetz and Albert (cf.
e.g. \cite{scor16},
\cite{lef}, \cite{albe}, \cite{ros}) and today one can find an exposition in
Chapters 5 and 9 of the book by Lange and Birkenhake \cite{l-b}, cf. also, for an
historical account, the article by Auslander and Tolimieri
\cite{a-t}.

Although the methods of Scorza and Rosati were more geometrrical, certainly more
than the later ones by Albert, who essentially worked in the new direction set up
by Emmy Noether, i.e., of the abstract algebra, a central role is played  by a
notion due to Rosati, the so-called {\bf Rosati involution}.

Given an
endomorphism with integral matrix $B$, the Rosati involution associates to it 
( $A$ being a Riemann integral matrix as in (3.3), (3.4)) the matrix $B' : A^{-1} \
^{t}B   A$. The Rosati involution is positive in the sense that the symmetric bilinear
form
$(B_1, B_2) := Tr (B_1' B_2 + B_1 B_2')$ yields a positive definite scalar product.

It turns out that the classification of Riemann matrices is very close to the
study of rational Algebras with a positive involution, and abstract arguments imply 
 that these  simple algebras are skew fields $\F$ of finite dimension over $\Q$
 of two types
 
 \begin{itemize}
 \item
 (I) The centre $K$ of $\F$ is a totally real number field and, if $K \neq \F$,
then $\F$ is a quaternion algebra over $K$. Moreover, for every embedding $\sigma:
K \ra \R$, $\F \otimes \R$ is always definite ($\F \otimes \R \cong \HH$), or
always indefinite ($ \F \otimes \R \cong M(2, \R)$).
\item
(II) The centre $K$ of $\F$ is a totally complex quadratic extension of  a totally 
real number field $K_0$, and then $\F \otimes \C$ is a matrix algebra $M(r,\C)$
such that the positive involution extends to the standard involution $ C \ra  \
\overline{ ^tC }$.
 \end{itemize}
 
 The analytic moduli theory of Abelian varieties owes much to the work of Siegel
and to his 'Symplectic geometry' (\cite{sieg43}): today the space of matrices
$\HHH_g : = \{\tau \in Mat(g,g,\C)| \tau = ^t\tau , \ Im (\tau) > 0 \} $ is called
the Siegel upper half space, and it is a natural parameter space for Abelian
varieties, since, depending on the polarization, there is a subroup $\Gamma$ of
$Sp(2g, \Q)$ such that the moduli space is, analytically, the quotient $ \HHH_g /
\Gamma$.

 The moduli theory of Abelian varieties with a certain polarization and
endomorphism structure was pursued relatively recently by Shimura (cf.
\cite{shim}), and it is a currently very active field of research for the arithmetic
applications of the theory of such Shimura varieties.

I do not need to cite for instance the
(recently proved) so called Shimura-Taniyama-Weil conjecture
about the modularity of elliptic curves defined over $\Q$: I can simply refer to
the talk by Wiles.

In this direction, however, the current tendency is to develop also much the
geometry,  since one has to look at the
reduction of these modular varieties modulo primes. The hope is that this study
will play a primary role for the pursuing of the so called Langlands program,
which is a vast generalization of the previously cited conjecture, proposing to
relate modular forms arising in different contexts ( cf. \cite{langl70}, and
\cite{langl76}, \cite{langl79}, \cite{del79} for early accounts of the story).

I 
hope that some more competent author than me will report about this
development in the present volume.

I want instead to end this section by pointing out (cf. \cite{zap}) how important
 the role of Scorza was for the development of the field of abstract algebra in
Italy: his path started with correspondences between curves, but, as we contended
here, his researches centered about Riemann matrices made him realize about the
relevance of the powerful new algebraic concepts.

\bigskip

\section{Further links to the italian school}

We mentioned in the previous section how the research of Rosati and Scorza was
very much influenced by the new geometric methods of the italian school of
algebraic geometry.

As we said, a crucial role was played by Castelnuovo: concerning Abel's theorem, in
the article
\cite{cast93}, entitled "The 1-1
correspondences between groups of $p$ points on
a curve of genus $p$", he explained how one could e.g. formulate the fundamental
theorem about the inversion of Abelian integrals as a consequence of the theory of
linear series on a curve (a development starting with the 
geometric
interpretation of the Riemann Roch theorem).

It is interesting to observe that, when he wrote a final note in the edition
appearing in his collected works, he points out that the results can be formulated
in a simpler way if one introduces the concept of the Jacobian variety of the
curve.

These notes added around 1935 are rather interesting: for instance, in the note to
the paper entitled "On simple integrals belonging to an irregular surface"
(\cite{cast05}) he pointed out that exactly in this Memoir he introduced the
concept of the so called {\bf Picard variety}, applyimg this concept to the
study of algebraic surfaces. In fact, the theorem of Picard to which Castelnuovo
refers, proved by Picard in \cite{pic}, and with precisions by Painlev\'e in
\cite{painl}, is indeed the characterization of the Abelian varieties, (we add to
it a slight rewording in modern language)

\begin{teo}
Let $V^p$ be a $p$-dimensional algebraic variety admitting a transitive
$p$-dimensional abelian group of birational transformations: then the points 
$\xi \epsilon V$ are uniformized by entire $2p$-periodic functions on $\C^p$, 
$$ \xi_k = \phi_k (u_1, \dots ,.. u_p) $$
(i.e., $V$ is birational to a complex torus of dimension $p$).
\end{teo}

The first main result of Castelnuovo in \cite{cast05} is  the construction of the so
called Albanese variety and Albanese map of an algebraic surface $X$.
Recall that, in modern language,  the Albanese variety of a projective variety $X$
is the Abelian variety
$ ( H^0 (\Omega^1_X ))^{\vee} / H_1(X)$ where $ H_1(X)$ is the lattice, in the dual
vector space of $ H^0 (\Omega^1_X )$, given by integration along closed paths.

The Albanese map, defined up to translation, as a result of the choice of a base
point $x_0$, associates to a point $x$ the linear functional $\alpha(x):=
\int_{x_0}^x  ( {\rm mod } \ H_1(X) )$.\footnote{ Giacomo Albanese emigrated from
Italy in 1936 to Sao Paulo, Brasil, where, soon after the war, he became closely
acquainted with Andr\'e Weil, who taught there, as well as Zariski. Weil is
responsible for the name "Albanese variety", but Ciliberto and Sernesi in
\cite{albcp} write: "Whilst the attribution of this concept to Albanese is
dubious,
..." . Indeed the basis for this is an article, \cite{alb34}, where
Albanese studies correspondences beween algebraic surfaces through the
consideration of the induced action on the space of holomorphic $1$-forms.
The coupled names "Albanese" and "Picard" appear in the title of the article by
W.L. Chow " Abstract theory of the Picard and Albanese varieties",\cite{chow}.

Indeed, in the 50's, one main  purpose was to distinguish among the two dual
varieties, which are only isogenous, and not in general isomorphic. According to
the historical note of Lang on page 52 of \cite{langav}, Matsusaka was the first
to give a construction of the Albanese variety using the generic curve.  }

The second result, obtained independently by F. Severi in \cite{sev05}, concerns
the equality of the irregularity of an algebraic surface and the dimension of the
space of holomorphic one forms. Both proofs were relying on a shaky proof given by
Enriques one year before, in \cite{enr05}, claiming the existence a continuous
system $\Sigma$  of dimension $q: = p_g - p_a$ of "inequivalent curves"(i.e., such
that for a generic curve $C \in \Sigma$, the set of curves in $\Sigma$ which are
linearly equivalent to $C$ has dimension zero).

Fortunately, a corrrect analytical proof was later found by Poincar\'e in
\cite{poi10}.  Enriques and Severi tried for a long time to repair the
flaw in Enriques' geometrical arguments, although in the end it started to
become clear the need for higher order differential elements (i.e., higher order
terms in the Taylor expansion of the curve variation). The fruit of the researches
carried on much later in the 50's was to show that indeed, for varieties defined
over algebraically closed fields of positive characteristic, the {\bf arithmetic
irregularity}
$q: = p_g - p_a$ was in general larger than the {\bf geometric irregularity}
defined as the maximal dimension of a continuous system of "inequivalent curves".
We refer the reader to Mumford's classical book
\cite{mumlc}, relating this question to the non reducedness of the {\bf Picard
scheme } $H^1( \hol_X^{*})$.

Castelnuovo was instead more interested in the applications of the previous
theorem, the most important one being the theorem \cite{cast05-2} that an
algebraic surface with  arithmetic genus $p_a$ smaller than $-1$ is birationally ruled.

This theorem is indeed one of the key theorems of the classification of algebraic
surfaces, since it also implies the well known 

{\bf Castelnuovo rationality
criterion: a surface is rational if and only if the bigenus $P_2$ and the
irregularity $q$ vanish}.

Without opening a new story, I would like to observe that the so called "Enriques
classification" of algebraic surfaces, done by Castelnuovo and Enriques, was one
of the most interesting cooperations in the history of mathematics, which took
place in the years from 1892 to 1914 (and especially intense in the period 1892 to
1906).

Besides the published papers, one can consult today the book entitled "Riposte
armonie" \cite{RA} ("Hidden harmonies", as are the ones governing algebraic
surfaces), which contains around 670 letters (or postcards) written from Federigo
Enriques to Guido Castelnuovo.

Naturally, also Castelnuovo wrote quite many letters, but apparently Enriques did
not bother to keep them. This already shows where  Castelnuovo and Enriques
respectively belong, in the rough distinction made by H. Weyl which we quoted
above, although both of them were obviously geometers. 

In fact, Enriques used to discuss mathematics with his assistents during long
walks in gardens or parks, and would only sometime stop to write something with
his stock on the gravel. Moreover, as Guido Castelnuovo wrote 
of him\footnote{In the
commemoration opening vol. I of his selected papers in geometry, \cite{enrMS},
IX - XXII.  }with affection , he was a "mediocre reader, who saw in a page
not what was written, but what he wanted to see"; 
certainly his brain was
always active like a vulcano. After the first world war the collaboration 
of the two broke
up, more on the side of Castelnuovo. As his daughter Emma Castelnuovo writes,
Enriques would regularly visit his sister 
(Castelnuovo's wife) at their house,
and after dinner the two mathematicians went 
to a separate sitting room, where
Enriques wanted to discuss his many new ideas, while Castelnuovo had prudently
instructed his wife to come
 after some time and interrupt their conversation with some excuse.
 
 Castelnuovo was 6 years older than Enriques, was always calm and mature, and,
after the appearance of his ground breaking two notes over algebraic surfaces
(\cite{cast91}) he was a  natural reference for the brilliant student Enriques,
who graduated in Pisa in 1891 (just at the time when Guido became a full
professor). Enriques wanted first to perfection his studies under the guidance
of Segre in Torino, but instead got a fellowship in Rome by Cremona, and there,
in 1892, started the intense mathematical interchange with Castelnuovo.  

In his first Memoir   (\cite{enr93}) Enriques, after an interesting historical
introduction, sketches the main tools to be used for the birational study of
algebraic surfaces, namely: the theory of linear systems of curves, the
canonical divisor and the operation of adjunction. Some results, as the claim
that two birational quartic surfaces are necessarily projectively equivalent, are
today known not to hold. It took a long time to make things work properly, and it
is commonly agreed that the joint paper \cite{c-e14} marks the achievement of the
classification theory.

For the later steps, (main ones, as well known to anyone who understands
the structure of the classification theorem) a very important role played the

{\bf(IP) De Franchis' theorem on irrational pencils }(\cite{df05},
\cite{cast05-2})  and the

{\bf(HS) Classification of hyperelliptic surfaces }.

{\bf(IP)} De Franchis' theorem, obtained independently also by Castelnuovo and
Enriques, asserts that if on a surface $S$ there are several linearly independent
one forms
$\omega_i, i=1, \dots r$ which are pointwise proportional, then there is a mapping
$f : S \ra C$ to an algebraic curve $C$ such that these forms are pull backs of
holomorphic one forms on $C$. It was used by Castelnuovo to show that under the
inequality $ p_g \leq 2q-4$ there is a mapping to a curve of positive genus.

This theorem leads to a typical example of the algebraization of the geometry:
although the Hodge theory of K\"ahler varieties was established with the use of
hard analytic tools which underlie Hodge's theorem on harmonic integrals
(\cite{hod}, \cite{kod}), what turned out to be very fertile was the simple
algebraic formulation in terms of the cohomology algebra of a projective variety.

Using this, Z. Ran, M. Green and myself (\cite{ran}, \cite{cat91}) were
independently able  to extend the result of de Franchis to the case of higher
dimensional varieties and higher dimensional targets. In this way the ideas of
the italian school came back to intense life, and became  an important tool for
the investigation of the fundamental groups of algebraic varieties (for instance
N. Mok \cite{mok} tried to extend this result to infinite dimensional
representations, with the hope of using such a result for the solution of the so
called Shafarevich conjecture about the universal coverings of algebraic
varieties).

{\bf(HS)} The second main work of de Franchis \footnote{Michele de Franchis, born
in Palermo, 1875, was very much influenced by the teaching of Scorza. He was
also the Director of the Rendiconti Circolo Matematico di Palermo in the years
1914-1946, as the successor of the founder G.B. Guccia }, together with Bagnera,
was the classification of hyperelliptic surfaces, i.e. of surfaces whose universal
covering space is biholomorphic to
$\C^2$ (cf.
\cite{ba-df07},
\cite{ba-df07-P} ,\cite{ba-df07}, \cite{e-s07},
\cite{e-s08},\cite{e-s09},\cite{e-s10}). This classification was also obtained by
Enriques and Severi, and the Bordin Prize was awarded to Enriques and Severi in
1907, and to the sicilian couple in 1909
\footnote{Lefschetz was instead recipient of the Bordin Prize in 1919 for the
cited Memoir \cite{lef}.}. Strange as it may seem that two couples get two prices
for the same theorem, instead of sharing one, this story is even more complicated,
since the first version of the paper by Enriques and Severi was withdrawn after a
conversation of Severi with de Franchis, and soon corrected. Bagnera and de
Franchis were only a little later, since they had to admit a restriction (a
posteriori useless, since no curve on an Abelian surface is contractible, cf.
\cite{df36}): their proof was however simpler, and further simplified by de
Franchis much later (\cite{df36-2}).

Another beautiful result, and this one even more related to Abelian integrals, is
the famous

{\bf Torelli' s Theorem: Let $C$, $C'$ be two algebraic curves whose Jacobian
varieties are isomorphic as polarized Abelian varieties (equivalently, admitting the
same matrix of periods for Abelian integrals of the I kind): then $C$ and $C'$ are
birationally isomorphic,
\cite{tor}}. 

Torelli was born in 1884 and was a student of Bertini, in Pisa, where he
attended also Bianchi's lectures: he was for short time assistant of Severi, and
died prematurely in the first world war, in 1915.

These years at the beginning of the 20-th century in Italy were thus quite
exciting. In the book by J. Dieudonn\'e, \cite{dieu} vol.1, Chapter VI, entitled
"Developpement et chaos", contains a paragraph dedicated to " L' \'ecole
italienne et la th\'eorie des systemes lin\'eaires", namely, devoted to the second
generation of the italian school. As third generation of the italian school, we
became very interested, also because of this criticism, to become fully
acquainted with the results of these precursors. 

What I found as a very interesting peculiarity reading the book by Enriques on
"Algebraic surfaces" (\cite{enr49}), was the mixture of theorems, proofs, speculations,
and history of the genesis of the mathematical ideas.

For instance, chapter IX, entitled "Irregular surfaces and continuous systems of
disequivalent curves", has a section 6 on "History of the theory of continuous
systems" (pages 339-347). One can read there that in 1902 Francesco Severi, who
had just graduated in Torino, following the advice of his master C. Segre,
accepted a position in Bologna as assistant of Severi. Under the influence of
Enriques, Severi started the investigation of the particular surfaces which occur
as the symmetric square of a curve (cf. \cite{sev03}). According to Enriques,
this research lead him to consider the {\bf Problem of the base} for the class of
divisors modulo numerical equivalence.

We may also observe that these surfaces have a very special geometry, and,
although a general characterization has not yet been described in general, we
have some quite recent (cf. \cite{ccm}, \cite{h-p01}, \cite{pi01}) results

\begin{teo}
The symmetric squares $C^{(2)}$ of a curve of genus $3$ are the only irregular
surfaces of general type with $p_g \geq 3$ presenting the {\bf non-standard case }
for the bicanonical map, i.e., such that 
\begin{itemize}
\item
Their bicanonical map $\phi_{2K}$ is not birational onto its image
\item
$S$ does not contain any continuous system of curves of genus $2$.
\end{itemize}

Moreover, any algebraic surface with $p_g = q = 3$ is either such a symmetric
product (iff $K^2=6$ for its minimal model), or has $K^2=8$ and is the quotient
of a product $C_1 \times C_2$ of two curves of respective genera $2$ and $3$ by
an involution $i = i_1 \times i_2$ where $C_1 / i_1$ has genus $1$, while $i_2$
operates freely.
\end{teo}

The symmetric square $C^{(2)}$ of a curve of genus $2$ occurs in another
characterization of the {\bf non-standard case } given by Ciliberto and
Mendes-Lopes (\cite{c-m1}, \cite{c-m2}).

\begin{teo}
The double covers $S$ of a principal polarized Abelian surface $A$,
branched on a divisor algebraically equivalent to $ 2 \Theta$, are the only
 irregular surfaces of general type with $p_g =2$ presenting the {\bf
non-standard case } for the bicanonical map.
\end{teo}

In these theorems plays an essential role the continuous system of {\bf
paracanonical curves }, i.e., of those curves which are algebraically equivalent
to a canonical divisor. To this system is devoted section 8 of the cited Chapter
IX of \cite{enr49}, and there Enriques, after mentioning false attempts by Severi
and himself to determine the dimension of the paracanonical system $\{ K\}$,
analyses the concrete case of a surface $C^{(2)}$, with $C$ of genus $3$, in
order to conjecture that $ dim \{ K\} = p_g$. The assumptions conjectured by
Enriques were not yet the correct ones, but, under the assumptions that the
surface does not contain any irrational pencil of genus $\geq 2$, the conjecture
of Enriques was proved by Green and Lazarsfeld, via the so called "Generic
vanishing theorems" (cf. \cite{g-l87}, \cite{g-l91}).

I will not dwell further on this very interesting topic , referring the reader to
the survey, resp. historical, articles \cite{cat91-B} \cite{cil91-B}.

I should however point out that further developments are taking place in this
direction, following a seminal paper by Mukai (\cite{muk81}) who extended the
concept of Fourier transforms (don't forget that theta functions are particular
Fourier series!) to obtain an isomorphism between the derived category of coherent
sheaves on an Abelian variety $X$ and the one of its dual Abelian variety
$\hat{X} := Pic^0(X)$.

One specimen is the combination of Mukai's technique with the theory of generic
vanishing theorems by Green and Lazarsfeld to obtain limitations on the
singularities of divisors on an Abelian variety (cf. \cite{hac00}, where one can
also find references to previous work by Kollar, Ein and Lazarsfeld).

Speaking about links with the italian school I should not forget the beautiful
lectures I heard in Pisa from Aldo Andreotti on complex manifolds and on complex
tori. Through his work with F. Gherardelli (\cite{a-g76}), I got in touch
with a problem of transcendental nature which occupies a central place in Severi's
treatise on Quasi Abelian Varieties.

{\bf Quasi Abelian varieties }, in Severi' s terminology, are the Abelian complex
Lie groups which sit as Zariski open sets in a projective variety.

Whence, they are  quotients $ \C^n / \Gamma$ where $\Gamma$ is a discrete
subgroup of $\C^n$, thus of rank $r \leq 2n$, and the above algebraicity
property leads  again to the two Riemann bilinear relations:

\begin{itemize}
\item
I) There exists an alternating form $ A : \Gamma \times \Gamma \ra \Z$
such that $A$ is the imaginary part of an Hermitian form $H$ on $V$
\item
II) there exists such an $H$ which is positive definite
(in this case $H$ is not uniquely determined by $A$).

\end{itemize}

Andreotti and Gherardelli conjectured that {\bf  The Riemann bilinear relations hold
if and only if there is a meromorphic function of $\C^n$ with group of periods
equal to $\Gamma$ }. This conjecture was the first Ph.D. problem I gave, and
after some joint efforts, Capocasa and I were able to prove it in \cite{c-c}.

\bigskip

{\bf  HOMOLOGICAL ALGEBRA for Abelian and irregular Varieties}

As already mentioned, David Mumford's ground breaking articles (\cite{mum66-7})
set up the scope of laying out a completely algebraic theory of theta functions.
His attempt was not the only one, for instance Barsotti (cf.
\cite{bar1},\cite{bar2},\cite{bar3},\cite{bar4}, and
\cite{langav} for references about his early work) had another approach to Abelian
varieties, based on power series, Witt vectors and generalizations of them (Witt
covectors cf.
\cite{bar3}), and the so called "Prostapheresis formula" (cf. \cite{bar4}).

Discussing here the respective merits of both approaches would be difficult, but
at least I can say that, while Barsotti's work is mainly devoted to Abelian
varieties in positive characteristic, the theory of theta groups, as already
mentioned, is also a very useful tool in characteristic zero.

The title of Mumford's series of articles is "On the equations defining abelian
varieties", which has a different meaning than "The equations defining abelian
varieties" \footnote{Computers however have not only helped us to do
some explicit calculations, but, according to F.O. Schreyer, they have also made us
wiser: by showing us explicit equations which  need several pages to be written down,
they make us wonder whether  having explicit equations  means any better
understanding.}. Thus, he set up a program which has been successfully carried out
in the case of several types of Abelian varieties. More generally, one can set as
a general target the one of studying the equations of irregular varieties, i.e.,
of those which admit a non trivial  Albanese map.

Since otherwise the problem is set in too high a generality, let me give a
concrete example (for many topics I will treat now, consult also  the survey
paper \cite{cat97}, which covers developments until 1996).

Chapter VIII of the book by F. Enriques (finished by his assistants Pompilj and
Franchetta after the death of Enriques) is devoted to the attempt to find
explicitly the canonical surfaces $S$ of low canonical degree in $\PP^3$: i.e.,
one considers surfaces with $p_g=4$ and with birational canonical map $\phi:= \phi_K
: S \ra \PP^3$. If the canonical system has no base points, then we will have a
surface of degree $d = K^2 \geq 5$. The cases $K^2 = 5,6$ are easy to describe,
and for $K^2 \geq 7$ Enriques made some proposals to construct some
regular surfaces (with a different method, Ciliberto \cite{cil81} was able to
construct these for $ 7 \leq K^2 \leq 10$, and to sketch a classification
program, later developed in \cite{cat84}, based essentially on Hilbert 's syzygy
theorem).

It was possible to treat the irregular case (cf. \cite{c-s02}) using a new
approach based on Beilinson's theorem (\cite{bei78}) for coherent sheaves $\FF$ on
$\PP^n$, which allows to write every such sheaf as the cohomology of a {\bf monad} (a
complex with cohomology concentrated at only one point) functorially associated
with
$\FF$.

The natural environment for rings not necessarily generated in degree $1$ is
however the {\bf weighted projective space}, which is the projective spectrum
$Proj(\A)$ of a polynomial ring $\A := \C[x_0, \dots, x_n]$ graded in a non
standard way, so that the indeterminates
$x_i$ have  respective degrees $m_i$ which are positive integers, possibly
distinct.

Canonaco (\cite{can00}) was first able to extend Beilinson's theorem to the
weighted case under some restriction on the characteristic of the base field,
and later (\cite{can02}) not only removed this restriction, but succeeded to
construct a unique functorial Beilinson type complex, making use of a new theory,
of so called {\bf graded schemes}. A concrete application given was to determine the
canonical ring of surfaces with $K^2=4, p_g = q = 2$ (for these surfaces there is
no good canonical map to the ordinary 3-space $\PP^3$).

Although a general theory appears to be very complicated, it thus turns  out that
Abelian varieties (and for instance related symmetric products of curves) offer
crucial examples (admitting sometime a geometrical characterization) for the
study of irregular varieties. Some of them were already discussed before, we want
to give a new one which is particularly interesting, and yields (cf. \cite{c-s02}
for more details) an easy counterexample to an old "conjecture" by Babbage
(cf. \cite{cat81} and \cite{bea79} for references and the first counterexample).

\begin{ex}
Let $J = J(C)$ be the Jacobian of a curve $C$ of genus $3$, and let $A \ra J$ be
an isogeny of degree $2$. The inverse image $S$ of the theta divisor $ \Theta
\cong_{bir} C^{(2)}$ provides $A$ with a polarization of type $(1,1,2)$.

The canonical map of $S$ factors through an involution $\iota$ with $32$ isolated
fixed points, and the canonical map of the quotient surface $\Sigma: = S /
\iota$, whose image is a surface  $Y$ of degree $6$ in $\PP^3$ having $32$ nodal
isolated singularities, and  a plane cubic as double curve. Moreover, $\Sigma$ is the
normalization of $Y$.
\end{ex}

It is now difficult for me to explain and to foresee exactly what principles
should these examples illustrate, let me however try. 

In Theorem 3.1 I tried to give an explicit example of the  "algebraization" of
the geometry, showing how the question of the birationality to Abelian varieties
of symmetric products of varieties can be reduced to pure exterior algebra
arguments. 
Abelian varieties are just, so to say, the complex incarnation of
exterior  algebras. 

On the other hand,  a companion article by
Bernstein-Gelfand-Gelfand (\cite{bgg78}) appeared next to the cited article
\cite{bei78} by Beilinson. In abstract setting, it shows the equivalence of the
derived category of coherent sheaves on $\PP^n$ and the category of finite
modules over the exterior algebras. While it is not yet completely clear how to
extend this result to the weighted case, quite recently in \cite{efs01} it was
shown how the BGG method allows to write functorially not only the sheaves, but
also the homomorphisms in the Beilinson monad. 

Note that the exterior algebra of
BGG is the exterior algebra over the indeterminates of a polynomial ring, and is
therefore apparently geometrically unrelated for the moment to the exterior
algebra of an Abelian variety. Given however a morphism $f : X \ra A$, we attach
to it the induced homomorphism between the respective holomorphic algebras of $A$
and $X$, $f^* : H^0(\Omega^*_A) \ra  H^0(\Omega^*_X)$. Whence, we obtain a
module over the exterior algebra, and we associate to it a Beilinson monad.
This procedere shows that to $f$ we associate some geometric objects related to
the Gauss maps corresponding to $f$.  

It is just the converse which looks more problematic, is it possible to associate,
to a map $f : X \ra \PP^n$ to a projective space, a geometric map to an Abelian
variety giving a realization of the module $f_* \hol_X$?  

In general, progress on the question of canonical rings or equations of
irregular varieties  might require at least further  combinations of  
the several existing  techniques which we have mentioned.

\bigskip

\section{More new results and open problems}

 {\bf 6.1. The Torelli problem}

The Torelli theorem, mentioned in the previous section, was again at the centre of
attention in the 50's, when several new proofs were found, by Weil, Matsusaka,
Andreotti ( and many others afterwards).

Particularly geometrical was the proof given by A. Andreotti (\cite{andr58}), who
showed that, given a Jacobian variety $J(C)$ with a theta divisor $ \Theta
\cong_{bir} C^{(g-1)}$, the Gauss map of $\Theta$ is ramified on the dual variety
$C^{\vee} \subset \PP^{g-1}$ of the canonical image of the curve $C$.
Using projective duality, one sees therefore that $(J(C), \Theta)$ determines $C$ 
(Andreotti had then  to treat  the hyperelliptic case separately).

Andreotti and Mayer (\cite{a-m67}) pushed the study of the geometry of canonical
curves, especially of the quadrics of rank $4$ containing them, to obtain some
equations valid for the period matrices of curves inside the Siegel upper half
plane. This paper, written at a time when the fashion was oriented in quite
different directions, had a great impact on the revival of classical researches
about algebraic curves, Abelian varieties (cf. \cite{acgh84} and
\cite{l-b} for references).\footnote{ I came to read this beautiful paper under
the instigation of Francesco Gherardelli. He explained to me  that the citation
from "Il teatro alla moda" by Benedetto Marcello: "As a first duty shall the
modern poet ignore all about the ancient Roman and Greek poets, because these
last too ignored everything about modern poetry" was motivated by the extreme
difficulties that Aldo had when he wanted to
 lecture in Paris about this result without first explaining the excellence of the
rings of coefficients he was using, or use the notation $g^1_r$ without explaining the
representability of such functors.}

In the same years Philip Griffiths (cf. \cite{grif68}, \cite{grif70},
\cite{grif84} and references therein) greatly extended the theory of the periods
of Abelian integrals, proposing to use the Hodge structures of varieties, i.e. the
isomorphism class of the datum
$$ H^*(X, \Z) \ra  H^*(X, \Z) \otimes \C = \oplus_{p,q} H^{p,q}(X)$$ to
study their moduli (note that a modern formulation of Torelli's theorem is that
the birational isomorphism class of an algebraic curve is determined by the Hodge
structure on its cohomology algebra).

A prominent role played in his program the attempt to find a reasonable
generalization of Torelli's theorem, and indeed (cf. \cite{grif84}) a lot of
Torelli type theorems were proved for very many classes of special varieties. 

Since then, a basic question has been the one of finding sufficient conditions for
the validity of an infinitesimal Torelli theorem for the period map of
holomorphic $n$-forms (i.e., for the Hodge structure on $ H^n(X, \Z) \otimes \C$)
for a variety $X$ of dimension $n$.  The question is then, roughly
speaking, whether the period map is a local embedding of the local moduli
space. Observe moreover that the
$n$-forms are the only forms which surely exist on simple cyclic ramified
coverings $Y_D$ of a variety $X$ of general type, branched on pluricanonical
divisors $D$ (divisors  $ D \in | mK_X |$), and in this context the Torelli
problem is a quantitative question about how large
$m$ has to be in order that the variation of Hodge structure distinguishes the
$Y_D$'s (cf. the article \cite{mig95} which gave  a very interesting application
of these ideas to families of higher dimensional varieties, opening a new
direction of research).

The validity of such an infinitesimal Torelli theorem can be formulated in purely
algebraic terms as follows: is the cup product 
$$ H^1(X, T_X) \times H^0 (X, \Omega^n_X) \ra H^1 (X, \Omega^{n-1}_X) $$ non
degenerate in the first factor?

Classical examples by Godeaux and Campedelli and modern ones
(cf.\cite{bpv84} for references) produce surfaces of general type with
$p_g = q=0$ (thus with a trivial Hodge structure),
yet with moduli.

In view of Andreotti's interpretation of Torelli's theorem, one suspects then that
a good condition might be the very ampleness of the canonical divisor, i.e., the
condition that the canonical map be an embedding.

This might unfortunately not be the case, as shown in a joint paper with I. Bauer
\cite{b-c02}:

\begin{teo}
There are surfaces of general type with onjective canonical morphism and 
such that 
the infinitesimal Torelli theorem for holomorphic $2$-forms does not hold for each
surface in the moduli space. Examples of such behaviour are quotients $(C_1 \times
C_2) / (\Z/3)$, where
$C_1 \ra C_1 / (\Z/3) \cong \PP^1$ is branched on $3k+2$ points, and $(\Z/3)$
acts freely on $C_2$ with a genus $3$ quotient.
\end{teo}

It would be interesting to establish stronger geometrical properties of the
canonical map which guarantee the validity of the infinitesimal Torelli theorem
for the 
holomorphic $n$-forms.

\bigskip

 {\bf 6.2. The Q.E.D. problem}
 
 If higher dimensional varieties were products of curves, life would be much
simpler. It obviously cannot be so, since there are plenty of varieties which are
simply connected (e.g., smooth hypersurfaces in $\PP^n$ with $n \geq 3$), without
being rational.

Can life be simpler ?

It is a general fact of life that, in order to make the study of algebraic
varieties possible, one must introduce some equivalence relation.

The most classical one is the so called birational equivalence, which allows   in
particular not to distinguish between the different projective embeddings of the same
variety.

\begin{df}
Let $X$ and $Y$ be projective varieties defined over the field $K$: then they are
said to be birational if their fields of rational functions are isomorphic : 
$K(X) \cong K(Y)$.
\end{df}

Moreover, one must allow algebraic varieties to depend on parameters, for instance
the  complex hypersurfaces of degree $d$ in $\PP^n$ depend on the coefficients of
their equations: but if these are complex numbers, we can have uncountably 
many birational classes of algebraic varieties.

To overcome this difficulty, Kodaira and Spencer introduced the notion of
{\bf deformation equivalence} for complex manifolds:
they
(\cite {k-s58}) defined two complex manifolds $X'$, $X$ to be
{\bf directly deformation equivalent} if there is a  proper holomorphic submersion
$ \pi : \Xi \ra
\Delta$ of a complex manifold $\Xi$ to the unit disk in the complex plane, such that
$X, X'$ occur as fibres of $\pi$. If we take the equivalence relation generated
by direct deformation equivalence, we obtain the relation of  deformation
equivalence, and we say that $X$ is a deformation of $X'$ in the large if
$X, X'$ are deformation equivalent.

These two notions extend to the case of compact complex manifolds the classical
notions of irreducible, resp. connected, components of moduli spaces. 

It was recently shown however( \cite{man01}, \cite{k-k01},\cite{cat01}, ) that it
is not possible to give effective conditions in order to guarantee the deformation
equivalence of algebraic varieties, as soon as the complex  dimension becomes
$\geq 2$.

Thus in \cite{catqed} I introduce the following relation

\begin{df}
Let $n$ be a positive integer, and consider, for complex algebraic varieties $X,
Y$ of dimension $n$, the equivalence relation generated by
\begin{itemize}
\item
(1) Birational equivalence
\item
(2) Flat deformations with connected base and with fibres having only at most
canonical singularities
\item
(3) Quasi \'etale maps, i.e., morphisms which are unramified in codimension $1$.

\end{itemize}
This equivalence will be denoted by $X \cong_{QED} Y$ (QED standing for:
quasi-\'etale-deformation). 
\end{df}

\begin{oss}

\begin{itemize}
\item
Singularities play here an essential role. Note first of all that, without the
restriction on these given in (2), we obtain the trivial equivalence relation
(since every variety is birational to a hypersurface).
\item
Assume that a variety $X$ is rigid, smooth with trivial algebraic fundamental
group: then $X$ has no deformations, and there is no non trivial quasi-\'etale map
$Y \ra X$. 

In this case the only possibility, to avoid that $X$ be isolated in its
QED-equivalence class, is that there exists a quasi-\'etale map $f : X \ra Y$.

If $f$ is not birational, however, the Galois closure of $f$ yields another
quasi-\'etale map
$\phi : Z \ra X$, thus it follows that $f$ is Galois and we have a contradiction
if $Aut(X) = 
\{1\}$.

It does not look so easy to construct such a variety $X$.
\end{itemize}
\end{oss}

Are there invariants for this equivalence?

A recent theorem of Siu (\cite{siu02})
shows that the Kodaira dimension is invariant by QED equivalence.

It is an interesting question to determine the QED equivalence classes inside the
class of varieties with fixed dimension $n$, and with Kodaira dimension $k$.
For curves and special surfaces, there turns out to be only one class
(\cite{catqed}):

\begin{teo}
In the case $ n \leq 2, k \leq 1$ the following conditions are equivalent
\begin{itemize}
\item 
(i)
$X \cong_{QED} Y$  
\item
(ii) $dim X= dim Y= n$, $Kod(X) = Kod(Y)= k$
\end{itemize}

\end{teo}
The previous result uses heavily the Enriques classification of algebraic
surfaces. We can paraphrase the problem, for the open case of surfaces of general
type, using Enriques' words (\cite{enr49}): "We used to say in the beginning
that, while curves have been created by God, surfaces are the work of the devil.
It appears instead that God wanted to create for surfaces a finer order of
hidden harmonies.."

Are here then  new hidden harmonies to be found?

\bigskip

{\bf Acknowledgements.}

I wish to thank  Arnfinn Laudal for convincing me to try to learn about Abel's
life and to look at his papers, and for allowing me to do it comfortably 
at home.

I would also like to thank Ragni Piene for some useful comments.

\bigskip

{\bf Bibliographical remark.}

We are not in the position to even mention the most important references.
However, some of the references we cite here contain a vast bibliography, for
instance \cite{sieg}, pages 193-240, and \cite{zar}, pages 248-268.

\vfill

\noindent
{\bf Author's address:}

\bigskip

\noindent 
Prof. Fabrizio Catanese\\
Lehrstuhl Mathematik VIII\\
Universit\"at Bayreuth\\
 D-95440, BAYREUTH, Germany

e-mail: Fabrizio.Catanese@uni-bayreuth.de


\begin{thebibliography}{99}


\bibitem[Abel81]{abel81}
N. H. Abel,
``{\em  Oeuvres compl\'etes, Tomes I et II }'', L. Sylow and S.Lie ed.,
 {\bf Im. Grondahl Son , Christiania } (1881), 621 pp., 341 pp. .
 
 \bibitem[Acc94]{acc}
 R.D.M. Accola, 
``{\em Topics in the theory of Riemann surfaces}'',
{\bf Springer Lecture Notes in Mathematics. 1595},   (1994),
ix, 105 p..
 
\bibitem[Alb34]{alb34}
G. Albanese,
``{\em  Corrispondenze algebriche fra i punti di due superficie algebriche. }'',
 {\bf Ann. R. Scuola Normale Sup. Pisa, (II), 3 } (1934), 1-26, 149-182 .


\bibitem[Alb-CP]{albcp}
G. Albanese,
``{\em  Collected papers. }'', with an historical note by C.Ciliberto and E.
Sernesi, {\bf Queen's papers in pure and applied mathematics, vol. 103 } (1996),
408 pp.


\bibitem[Albe]{albe}
A.A. Albert,
``{\em  Structure of algebras. }''
{ \bf  Amer. Math. Soc. Coll. Publ. 24, New York}, (1939).

\bibitem[Andr58]{andr58}
A. Andreotti, 
``{\em  On a theorem of Torelli. }''
{ \bf Am. J. Math. 80}, (1958), 801-828 .

\bibitem[A-M67]{a-m67}
A. Andreotti, A. Mayer,
``{\em  On period relations for abelian integrals on algebraic curves}''
{ \bf Ann. Scuola Normale Sup. Pisa 21}, (1967), 189-238.

\bibitem[A-G76]{a-g76}
A. Andreotti, F. Gherardelli,
``{\em  Some remarks on quasi-abelian manifolds}'',
in 'Global analysis and its applications',
{ \bf I.A.E.A. Vienna}, (1976), 203-206.

\bibitem[App-Gours95]{appgour}
P. Appell, E. Goursat
``{\em  Th\'eorie des fonctions alg\'ebriques et de leurs int\'egrales. }''
{ \bf Gauthier-Villars, Paris},(1895), 530 pp..

\bibitem[ACGH84]{acgh84}
E. Arbarello, M. Cornalba, P. Griffiths, J. Harris,
``{\em  Geometry of algebraic curves I. }''
{ \bf Grundlehren der math. Wiss. 267, Springer},(1984), 386 pp..

\bibitem[Ati57]{ati57}
 M. F. Atiyah, 
``{\em Vector bundles over an elliptic curve. }''
{ \bf Proc. Lond. Math. Soc., III. Ser. 7}, 414-452 (1957). 

\bibitem[A-T]{a-t}
L. Auslander, R. Tolimieri
``{\em A matrix free treatment of the problem of Riemann matrices. }''
{ \bf Bull. Amer. Math. Soc. 5, 3  }, (1981), 263-312.

\bibitem[Ba-DF07-P]{ba-df07-P}
L. Bagnera, M. de Franchis
``{\em Sur les surfaces hyperelliptiques. }''
{ \bf Comptes Rendus Acad. Sci. Paris, 145 }, (1907), 747-749.

\bibitem[Ba-DF07]{ba-df07}
L. Bagnera, M. de Franchis
``{\em Sopra le superficie algebriche che hanno le coordinate del punto generico
esprimibili con funzioni meromorfe quadruplamente periodiche di due parametri, I,
II. }'' {\bf Rend. Accad. Lincei, 16 }, (1907), 492-498, 596-603.

\bibitem[Ba-DF08]{ba-df08}
L. Bagnera, M. de Franchis
``{\em Le superficie algebriche le quali ammettono una rappresentazione
parametrica mediante funzioni iperellittiche di due argomenti. }'' {\bf
Memorie Soc. Italiana delle Scienze, detta Accad. dei XL, 15 }, (1908),
251-343.

\bibitem[Bars64-6]{bar1}
I. Barsotti,
``{\em Medoti analitici per variet\'a abeliane in caratteristica positiva. Cap. 1- 2 ,
Cap. 3-4 , Cap.  5, 6, 7 }'' {\bf
Ann. Sc. Norm. Super. Pisa, Sci. Fis. Mat., III. }{Ser. 18, } (1964),
1-25 , {Ser. 19, } (1965), 277-330, (1965),
  481- 512 ; {Ser. 20, }(1966), 101-137, 331-365 .
  
  \bibitem[Bars70]{bar2}
I. Barsotti,
``{\em Considerazioni sulle funzioni theta ,}'' 
{\bf Sympos. Math., Roma 3, Probl. Evolut. Sist. solare, 
Nov. 1968 e Geometria, Febb. 1969}, (1970), 247-277 .

 \bibitem[Bars81]{bar3}
I. Barsotti,
``{\em Bivettori,}'' in 'Algebraic geometry, int. Symp. Centen. Birth F. Severi',
Roma 1979, {\bf Symp. Math. 24}, (1981), 23-63 . 

 \bibitem[Bars83-5]{bar4}
I. Barsotti,
``{\em Differential equations of theta functions, Continuation.}''
{\bf Rend. Accad. Naz. Sci. Detta XL, V. Ser., Mem. Mat. 7, No.1.} , (1983), 227-276 
and  {\bf Mem. Mat. 9}, (1985), 215-236 . 





\bibitem[BPV84]{bpv84} 

W. Barth, C. Peters and A. Van de Ven,
``{\em Compact Complex Surfaces}'', {\bf Ergebnisse der
Mathematik und ihrer Grenzgebiete, Folge 3, B.4, Springer-Verlag}, (1984), 304
pp..

\bibitem[B-C02]{b-c02}
I. Bauer, F. Catanese,
``{\em  Symmetry and variation of Hodge structure}'', preprint 2002.

\bibitem[Bea79]{bea79}  
A. Beauville, 
``{\em L'application canonique pour les
                 surfaces de type general}'',{\bf  Inv. Math. 55 },(1979),
                 121-140.

\bibitem[Bei78]{bei78}
A. Beilinson,
``{\em  Coherent sheaves on $\PP^n$ and problems of linear algebra}'',
{ \bf Funkt. Anal. i Pril. 12 (3)}, (1978), 68-69, translated in {\bf Funct.
Anal. Appl. 12 }, 214-216.

\bibitem[BGG78]{bgg78}
I.N. Bernstein, I.M. Gelfand, S.I. Gelfand,
``{\em  Algebraic bundles on $\PP^n$ and problems of linear algebra}'',
{ \bf Funkt. Anal. i Pril. 12 (3)}, (1978), 66-68, translated in {\bf Funct.
Anal. Appl. 12 }, 212-214.

\bibitem[Bia16]{bia}
L. Bianchi,
``{\em Lezioni sulla teoria delle funzioni di variabile complessa}'',
{ \bf Spoerri, Pisa}, (1916).

\bibitem[Bott02]{bott02}
U. Bottazzini, 
``{\em ``Algebraic truths'' vs ``geometric fantasies'': Weierstrass' response to
Riemann. }''  Li, Ta Tsien (ed.) et al., Proceedings of the international congress
of mathematicians, ICM 2002, Beijing, China, August 20-28, 2002. Vol. III: Invited
lectures.{\bf  Beijing: Higher Education Press}. 923-934 (2002).

\bibitem[Can00]{can00}
A.Canonaco,
``{\em  "A Beilinson-type theorem for coherent sheaves on weighted projective
spaces }'', { \bf Jour. of Algebra 225}, (2000), 28-46.


\bibitem[Can02]{can02}
A.Canonaco,
``{\em  "The Beilinson complex and canonical rings of irregular surfaces"}'',
{ \bf Tesi di Perfezionamento, Scuola Normale Superiore, Pisa}, (2002).

\bibitem[C-C91]{c-c}
F. Capocasa - F. Catanese, 
``{\em Periodic meromorphic functions.}''
{\bf Acta Math. 166} (1991) 27-68.

\bibitem[Cast91]{cast91}
G. Castelnuovo, 
``{\em Osservazioni intorno alla geometria sopra una superficie. Note I, II}''
{\bf Rendiconti del R. Ist. Lombardo, s.2, vol. 24} (1891) pages 246-254 and
255-265 of the Memorie Scelte.

\bibitem[Cast93]{cast93}
G. Castelnuovo, 
``{\em Le correspondenze univoche tra gruppi di $p$ punti sopra una curva di
genere $p$.}'' {\bf Rendiconti del R. Ist. Lombardo, s.2, vol. 25} (1893) pages
79-94 of the Memorie Scelte.

\bibitem[Cast05]{cast05}
G. Castelnuovo, 
``{\em Gli integrali semplici appartenenti ad una superficie irregolare.}'' {\bf
Rendiconti della R. Accad. dei Lincei, s.V, vol. XIV} (1905) pages 473-500 of the
Memorie Scelte.

\bibitem[Cast05-2]{cast05-2}
G. Castelnuovo, 
``{\em Sulle superficie aventi il genere aritmetico negativo.}'' {\bf
Rendiconti del Circ. Mat. di Palermo, t. XX } (1905) pages 501-506 of the
Memorie Scelte.


\bibitem[Cast]{cast}
G. Castelnuovo, 
``{\em Memorie scelte.}'' {\bf Zanichelli, Bologna}
(1937).

\bibitem[C-E14]{c-e14}
G. Castelnuovo, F. Enriques
``{\em Die algebraischen Fl\"achen vom Gesichtspunkte der birationalen
Transformationen aus.}'' {\bf  Enzycl. der math. Wissensch. III , 2, 6-6b}
(1914), 674-768.


\bibitem[Cat81]{cat81}
F. Catanese, 
``{\em Babbage's conjecture, contact of surfaces, symmetric determinantal
varieties and applications.}'', {\bf
Inv. Math. 63} (1981), 433-465. 

\bibitem[Cat84]{cat84}
F. Catanese, 
``{\em Commutative algebra methods and the equations of regular surfaces.}'' 
in 'Algebraic Geometry Bucharest 1982', {\bf Springer Lect. Notes in Math. 1056} (1981),
68-111. 


\bibitem[Cat91]{cat91}
F. Catanese, 
``{\em Moduli and classification of irregular K\"ahler manifolds
 (and algebraic varieties) with Albanese general type fibrations.
 Appendix by Arnaud Beauville.}''
{\bf Inv. Math. 104} (1991) 263-289; Appendix 289 .

\bibitem[Cat91-B]{cat91-B}
F. Catanese, 
``{\em Recent results on irregular surfaces and irregular K\"ahler manifolds.}''
in 'Geometry and complex variables'
{\bf Lecture Notes in pure appl. Math. 132,  Marcel Dekker, New York} (1991),
59-88.

\bibitem[Cat97]{cat97}
F. Catanese,
``{\em Homological algebra and algebraic surfaces}'',
in 'Algebraic Geometry', J.Koll\'ar. et al. ed. {\bf  Proc. Symp. Pure Math. 62,
Amer. Math. Soc., Providence R.I. } (1997) , 3 - 56.

\bibitem[Cat02]{cat02}
F. Catanese,
``{\em Deformation types of  real and complex manifolds}'',
in the Proc. of the Chen-Chow Memorial Conference, 
Advanced Topics in Algebraic Geometry and Algebraic
Topology, {\bf  Nankai Tracts in Math. 5, World
Scientific } (2002) ,  193-236.

\bibitem[Cat01]{cat01}
F. Catanese, 
``{\em  Moduli  Spaces of Surfaces and Real Structures}'', 
{\bf Ann. Math. 158 } (2003), n.2, 539-554.

\bibitem[Cat02-2]{cat02-2}
F. Catanese,
``{\em Deformation in the large of some complex manifolds, I}'',
preprint 2002, 
to appear in a Volume in Memory of Fabio Bardelli, Ann. Mat. pura e appl. 


\bibitem[CatQED]{catqed}
F. Catanese,
{\em Q.E.D. for algebraic varieties},
Preprint (2002). 

\bibitem[Ca-Ci93]{ca-ci}
F. Catanese, C. Ciliberto, 
``{\em Symmetric products of elliptic curves and surfaces of general type with $p_g =
q=1$.}'' { \bf J. Algebr. Geom. 2, No.3}, 389-411 (1993).

\bibitem[C-C-ML98]{ccm}
F. Catanese, C. Ciliberto, Mendez-Lopez,
``{\em On the classification of irregular surfaces of general type with
  nonbirational bicanonical map,}'' { \bf
Trans. Am. Math. Soc. 350, No.1},  (1998), 275-308.




\bibitem[C-S02]{c-s02}
F. Catanese, F.O. Schreyer
{\em Canonical projections of irregular algebraoc surfaces},
in 'Algebraic Geometry. A volume in memory of Paolo Francia,
{\bf De Gruyter, Berlin, New York }(2002), 79-116. 

\bibitem[Cay78-1]{cay1} 
A. Cayley, 
``{\em Desiderata and suggestions n.1. The theory of groups}'',
{\bf Amer. Jour. Math, Vol. 1, 1} (1878), 50-52.

\bibitem[Cay78-2]{cay2} 
A. Cayley, 
``{\em Desiderata and suggestions n. 2. The theory of groups: graphical
representation}'', {\bf Amer. Jour. Math, Vol. 1, 2} (1878), 174-176.

\bibitem[Chab50]{chab}
 C. Chabauty, 
``{\em Limite d'ensembles et gŽomŽtrie des nombres}'',{\bf Bull. Soc.
Math. France 78}, (1950), 143--151.


\bibitem[Che58]{che58} 
S.S. Chern, 
``{\em Complex manifolds}'',
{\bf Publ. Mat. Univ.  Recife} (1958).

\bibitem[Chow59]{chow} 
W.L. Chow, ``{\em Abstract theory of the Picard and Albanese varieties}'',
{\bf Ann. of Math. }, (1959),  .


\bibitem[Cil81]{cil81} 
C. Ciliberto, 
``{\em Canonical surfaces with $p_g=p_a=4$ and
     $K^2 = 5,\dots,10$,}'', {\bf  Duke Math. J. 48}, (1981), 121-157.


\bibitem[Cil91-B]{cil91-B}
C.Ciliberto, 
``{\em A few comments on some aspects of the mathematical work of F. Enriques.}''
in 'Geometry and complex variables'
{\bf Lecture Notes in pure appl. Math. 132,  Marcel Dekker, New York} (1991),
89-109.

\bibitem[Cil98]{cil98}
C. Ciliberto, 
``{\em M. de Franchis and the theory of hyperelliptic surfaces}'',
in 'Studies in the history of modern mathematics',
{\bf Supplemento Rend. Circ. Mat. di Palermo, S.II 55 }(1998), 45-73.


\bibitem[Ci-ML00]{c-m1}
C. Ciliberto, M.Mendes Lopes,
``{\em On surfaces with $p_g = q =2$ and non birational bicanonical map}'',
{\bf Adv. in Geom. 2} (2002), n. 3, 281-300.

\bibitem[Ci-ML02]{c-m2}
C. Ciliberto, M.Mendes Lopes,
``{\em On surfaces with $p_g = 2, q =1$ and non-birational bicanonical map}'',
in 'Algebraic Geometry. A volume in memory of Paolo Francia,
{\bf De Gruyter, Berlin, New York }(2002), 117-126. 

\bibitem[Ci-Se]{ci-se}
C. Ciliberto, E. Sernesi,
``{\em Some aspects of the scientific activity of Michele de Franchis}'',
in 'Opere di de Franchis',
{\bf Supplemento Rend. Circ. Mat. di Palermo, S.II 27 }(1991), 3-36.

\bibitem[Coble]{coble}
 Coble, Arthur B. 
``{\em  Algebraic geometry and theta functions}''. Revised printing.
 {\bf American
Mathematical Society Colloquium Publication, vol. X American Mathematical Society,
Providence, R.I}, (1929), (II edition 1961) vii+282 pp . 

\bibitem[corn76]{corn}
M. Cornalba, 
``{\em Complex tori and Jacobians}''.
{\bf Complex Anal. Appl., int. Summer Course Trieste 1975, Vol. II}, (1976) 39-100 .

\bibitem[Cous02]{cous02}
P.  Cousin, 
``{\em  Sur les fonctions periodiques}''
{\bf Ann. Sci. \'Ecole Norm. Sup. 19} (1902), 9-61.
 
\bibitem [Dav79-1]{dav79-1}

Davenport, J.H.
''{\em  Algorithms for the integration of algebraic functions.  }''
Symbolic and algebraic computation, EUROSAM '79, int. Symp., 
Marseille 1979,{\bf  Lect. Notes Comput. Sci. 72}, 415-425 (1979). 

\bibitem [Dav79-2]{dav79-2}

Davenport, J.H.
''{\em  The computerisation of algebraic geometry. }'' 
Symbolic and algebraic computation, EUROSAM '79, int. Symp.,
Marseille 1979,{\bf Lect. Notes Comput. Sci. 72}, 119-133 (1979).


\bibitem [Dav81]{dav81} 

Davenport, James Harold
''{\em On the integration of algebraic functions. }''
{\bf Lecture Notes in Computer Science, 102. Berlin Heidelberg New York:
Springer-Verlag} 197 p.  (1981).

\bibitem [Dini78]{dini}
U. Dini,
''{\em Fondamenti per la teorica delle funzioni di variabile reale. }''
{\bf Spoerri, Pisa}, (1878). 

\bibitem[D-F05]{df05}
M. de Franchis,
''{\em Sulle superficie algebriche le quali contengono un fascio irrazionale di
curve }'', {\bf Rendiconti Circ. Mat. di Palermo, XX}(1905), 49-54.

\bibitem[D-F36]{df36}
M. de Franchis,
''{\em Dimostrazione del teorema fondamentale sulle superficie iperellittiche }'',
{\bf Rendiconti Accad. Lincei, 24}(1936), 3-6.

\bibitem[D-F36-2]{df36-2}
M. de Franchis,
''{\em Sulla classificazione delle superficie iperellittiche }'',
{\bf Scritti matematice in onore di Berzolari}(1936), 613-615.


\bibitem[Dieu]{dieu}
J. Dieudonn\'e,
''{\em Cours de g\'eom\'etrie alg\'ebrique, vol. 1,2 }'', {\bf Collection SUP
10,11, Presses Universitaires de France }(1974), 234 pp. 222 pp.,  

\bibitem[Del79]{del79}
P. Deligne,
''{\em Vari\'et\'es de Shimura : interpr\'etation modulaire, et techniques de
construction de mod\'eles canoniques. }'', in 'Automorphic forms, representations,
and L-functions', {\bf Proc.Symp. Pure Math. 33, vol. 2, A.M.S.}(1979), 247-290.

\bibitem[Don81]{don}
R. Donagi, 
''{\em  The tetragonal construction}'',
{\bf Bull. Am. Math. Soc., New Ser. 4}, (1981) 181-185 .

\bibitem[E-F-S]{efs01} D. Eisenbud, G. Fl\o ystad, F.-O. Schreyer,
   ''{\em Sheaf cohomology and free fesolutions over exterior algebras}'',
    preprint,  math.AG/0104203.

\bibitem[EDM1]{edm1}

 ``{\em Encyclopedic dictionary of mathematics. Vol. I.
   Abel to multivariate analysis. }'' Translated from the second Japanese
edition. Edited by Sh™kichi Iyanaga and Yukiyoshi Kawada. {\bf MIT Press, Cambridge,
Mass.-London}, (1977). xv+883 pp. 

\bibitem[EDM2]{edm2}
``{\em Encyclopedic dictionary of mathematics. Vol. II.
Networks to zeta functions.}'' Translated from the second Japanese edition. Edited by
Sh™kichi Iyanaga and Yukiyosi Kawada. {\bf  MIT
Press,
   Cambridge, Mass.-London}, 1977. pp. 885--1750. 

\bibitem[Enr93]{enr93}

F. Enriques,
''{\em Ricerche di geometria sulle superficie algebriche.}'' 
{\bf   Memorie Accad. Torino, s.2 vol. 44}(1893), 171-232.

\bibitem[Enr05]{enr05}

F. Enriques,
''{\em Sulla proprieta' caratteristica delle superficie irregolari.}'' 
{\bf   Rend. Accad. Sci. Bologna, nuova Serie vol. 9}(1905), 5-13.

\bibitem[Enr49]{enr49}

F. Enriques,
''{\em Le superficie algebriche.}'' 
{\bf   Zanichelli, Bologna}, (1949), 464 pp. (Translation by F.Catanese, C.
Ciliberto, R. Pardini, to be published by Canbridge University Press).

\bibitem[EnrMS]{enrMS}

F. Enriques,
''{\em Memorie scelte di geometria, vol. I, II, III.}'' 
{\bf   Zanichelli, Bologna}, (1956), 541 pp.,(1959), 527 pp.,(1966), 456 pp. .

\bibitem[E-S07]{e-s07}

F. Enriques, F. Severi,
''{\em Intorno alle superficie iperellittiche.}'' 
{\bf   Rend. Accad. Lincei, s. V vol. 16}(1907), 443-453.

\bibitem[E-S08]{e-s08}

F. Enriques, F. Severi,
''{\em Intorno alle superficie iperellittiche irregolari.}'' 
{\bf   Rend. Accad. Lincei, s. V vol. 17}(1908), 4-9.

\bibitem[E-S09]{e-s09}

F. Enriques, F. Severi,
''{\em M\'emoire sur les surfaces hyperelliptiques, I.}'' 
{\bf   Acta Math.  vol. 32}(1909), 283-392.

\bibitem[E-S10]{e-s10}

F. Enriques, F. Severi,
''{\em M\'emoire sur les surfaces hyperelliptiques, II.}'' 
{\bf   Acta Math.  vol. 33}(1910), 223-403.

\bibitem[Falt83]{falt83}

G. Faltings,
''{\em EndlichkeitssŠtze fŸr abelsche VarietŠten Ÿber Zahlkšrpern.}'' 
{\bf   Invent. Math. 73 ,no. 3}(1983),  349--366.

''{\em Erratum: "Finiteness theorems for abelian varieties over number fields}''.
{\bf Invent. Math. 75, no. 2} (1984), 381.


 \bibitem[Falt-Ch]{falt-ch}
 
 G. Faltings, C.L.  Chai, ''{\em  Degeneration of abelian
varieties.}'' 
  {\bf  Ergebnisse der Mathematik und ihrer
Grenzgebiete (3) , 22.
Springer-Verlag, Berlin}, 
  ( 1990). xii+316 pp. 
  
  \bibitem [Fisch86]{fisch86}
 G. Fischer, (ed.)
''{\em  Mathematische Modelle.}''
Aus den Sammlungen von Universit\"aten und Museen. Mit 132 Fotografien.
(Mathematical models. From the collections of universities and museums. With 132
photographs). Bildband und Kommentarband. {\bf Braunschweig/Wiesbaden: Friedr.
Vieweg \& Sohn},
 (1986), XII, 129 pp.; VIII, 89 pp..
 
  \bibitem[G-P84]{g-p}
A. Genocchi, 
''{\em Calcolo differenziale}'', con aggiunte del prof. G. Peano {\bf Torino}, (1884).

 
 \bibitem[G-L87]{g-l87}
M. Green, R. Lazarsfeld 
''{\em Deformation theory, generic vanishing theorems, and some conjectures of
Enriques, Catanese and Beauville}''. {\bf Inv. Math. 90}, (1987), 416-440.

\bibitem[G-L91]{g-l91}
M. Green, R. Lazarsfeld 
''{\em Higher obstructions to deforming cohomology groups of line bundles}''. {\bf
Jour. Amer. Math. Soc. 4}, (1991), 87-103.

\bibitem[Grif68]{grif68}
P. Griffiths,  
''{\em Periods of integrals on algebraic manifolds I,II}''. {\bf
Amer. Jour. Math. 90}, (1968), 568-626, 805-865.

\bibitem[Grif70]{grif70}
P. Griffiths,  
''{\em Periods of integrals on algebraic manifolds III}''. {\bf
Publ. Math. I.H.E.S. 38}, (1970), 125-180.


\bibitem[Grif84]{grif84}
P. Griffiths, (ed) 
''{\em Topics in transcendental algebraic geometry}''. {\bf
Annals of Math. Studies 106, Princeton Univ. Press}, (1984), 316 pp..

\bibitem[Hac00]{hac00}
C.D. Hacon, 
''{\em Fourier transforms, generic vanishing theorems and polarizations of
abelian varieties}''. {\bf  Math. Zeit. 235}, (2000), 717-726.

\bibitem[H-P01]{h-p01}
C.D. Hacon, R. Pardini,
''{\em Surfaces with $p_g = q =3$}'',  {\bf  Trans. Amer. Math. Soc. 354} 
(2002), n.7, 2631-2638.

\bibitem[Hod52]{hod}
W.V.D. Hodge,
''{\em Theory and application of the harmonic integrals}'',  {\bf 
Cambridge Univ. Press} (1952).


\bibitem[Igu]{igu}
J.I. Igusa, 
''{\em Theta functions }''. {\bf Die Grundlehren der mathematischen
Wissenschaften, Band 194. Springer-Verlag, New York-Heidelberg}, (1972) x+232 pp.


  
\bibitem[Kempf73]{kempf73}  
 G.Kempf, ''{\em On the geometry of a theorem of Riemann}'' . 
 { \bf Ann. of Math. (2)
98 }(1973), 178--185. 

\bibitem[KempfRINGS]{kempfrings}  
G. Kempf,  ''{\em 
Projective coordinate rings of abelian varieties.}'' 
{ \bf Algebraic analysis, geometry, and number theory (Baltimore, MD, 1988)
  Johns Hopkins Univ. Press, Baltimore, M}, (1989), 225--235. 
  
  \bibitem[K-K01]{k-k01}  V. Kharlamov and V. Kulikov,''{\em 
On real structures of real surfaces . }'' 
{\bf math.AG/0101098}.

 \bibitem[Klein78]{klein78}
F. Klein,
''{\em 
\"Uber die Transformationen siebenter Ordnung der elliptischen Funktionen}'',
{\bf Math.Ann. 14 }, (1878)  .

 \bibitem[KleinCW]{klein}
 F. Klein, 
``{\em  Gesammelte mathematische Abhandlungen , B\"ande I, II, III}''.

    Band I: Liniengeometrie, Grundlegung der Geometrie, zum Erlanger Programm.
Edited by R. Fricke und A. Ostrowski (with additions by Klein), ii+xii+612 pp. (1
plate).

   Band II: Anschauliche Geometrie, Substitutionsgruppen und
Gleichungstheorie, zur mathematischen Physik. Edited by  R. Fricke und H.
Vermeil (with additions by Klein), ii+vii+713 pp. 

Band III : 
 Elliptische Funktionen, insbesondere Modulfunktionen, hyperelliptische und Abelsche
Funktionen, Riemannsche Funktionentheorie und automorphe
   Funktionen, Anhang verschiedene Verzeichnisse. Edited by
    R. Fricke, H. Vermeil und E. Bessel-Hagen (with additions by F. Klein),
ii+ix+774+36 pp.

 Reprint of the first edition {\bf Julius Springer Verlag, Berlin},
(1923), {\bf Springer-Verlag,
   Berlin-New York}, 1973.

\bibitem [Kod52]{kod}
K. Kodaira ,``{\em The theory of harmonic integrals and their applications to
algebraic geometry}'' (1952), in "Kunihiko Kodaira: Collected Works", {\bf
Princeton Univ. Press and Iwanami Shoten}, (1975), 488-582 .

\bibitem [K-M71]{k-m71}
K. Kodaira, J. Morrow ,``{\em Complex manifolds}'' {\bf Holt,
Rinehart and Winston},  New York-Montreal, Que.-London (1971).

\bibitem[K-S58]{k-s58}
K. Kodaira , D. Spencer ''{\em On deformations of complex analytic
structures I-II}'', {\bf Ann. of Math.  67} (1958), 328-466 .

\bibitem[Kraz]{kraz}
A. Krazer,
``{\em Lehrbuch der Thetafunktionen}'' 
{\bf  Teubner, Leipzig}, (1903).

\bibitem[K-W]{k-w}  
A. Krazer, W. Wirtinger
``{\em Abelsche Funktionen und allgemeine Thetafunktionen}'' 
{\bf Enzykl. Math. Wiss. II B 7}, (1920), 604-873.

\bibitem[Lang-AV]{langav}
S. Lang, C. 
''{\em  Abelian varieties. }''
{\bf Interscience Publ. Inc. }, (1959), reprinted by Springer Verlag (1983).


\bibitem[L-B]{l-b}
H. Lange, C. Birkenhake,
''{\em Complex Abelian varieties. }''
{\bf Grundlehren der mathematischen Wissenschaften 302, Springer Verlag }, (1992).

\bibitem[Langl70]{langl70}
R.P. Langlands,
''{\em Problems in the theory of automorphic forms. }'', in 'Lectures in modern
analysis and applications', {\bf Springer Lecture Notes in Math. 170}(1970), 18-86.

\bibitem[Langl76]{langl76}
R.P. Langlands,
''{\em Some contemporary problems with origins in the Jugendtraum.
}'', in 'Mathematical developments arising from Hilbert problems', {\bf Proc.Symp.
Pure Math. 28, vol. 2, A.M.S.}(1976), 401-418.

\bibitem[Langl79]{langl79}
R.P. Langlands,
''{\em Automorphic representations, Shimura varieties, and motives, ein M\"archen.
}'', in 'Automorphic forms, representations, and L-functions', {\bf Proc.Symp.
Pure Math. 33, vol. 2, A.M.S.}(1979), 205-246.

\bibitem[Leb02]{leb}
H. Lebesgue, 
''{\em Int\'egrale, longuer, aire. }'',
{ \bf Annali di Mat. pura ed appl. 7},(1902), 231-358.

\bibitem[Lef21]{lef}
S. Lefschetz, 
''{\em On certain numerical invariants of algebraic varieties with
application to abelian varieties }'',
{ \bf Trans. Amer. Math. Soc. 22 } no. 3,(1921),  327--406  and no. 4,
(1921),407--482. 

\bibitem[Man01]{man01}
M. Manetti, ''{\em On the Moduli Space of diffeomorphic algebraic
surfaces}'',
{\bf  Inv. Math. 143} (2001), 29-76. 


\bibitem[Man-Schr]{man-schr}
N. Manolache, F.O.  Schreyer, 
''{\em  Moduli of $(1,7)$-polarized abelian surfaces via syzygies. }'' 
{ \bf Math. Nachr. 226}, 177-203 (2001).

\bibitem[Matsu58]{matsu}
T. Matsusaka, 
''{\em On a theorem of Torelli. }''
{ \bf Am. J. Math. 80}, 784-800 (1958).

\bibitem[Mig95]{mig95}
L. Migliorini, 
''{\em A smooth family of minimal surfaces of general type over a curve of genus
at most one is trivial. }'' { \bf Jour. Alg. Geom. 4}, (1995), 668-684.

\bibitem[Mok97]{mok}
N. Mok, 
''{\em The generalized theorem of Castelnuovo-de Franchis for unitary
representations.}''  in '  Geometry from the Pacific
Rim. Proceedings of the Pacific Rim geometry conference, National University of
Singapore, Republic of Singapore, December 12--17, 1994,' Berrick, et al.,ed., {\bf  
Walter de Gruyter, New
York} (1997), 261-284.

\bibitem[Muk81]{muk81}
S. Mukai, 
''{\em Duality between $D(X) and D(\hat{X})$ with its applications to Picard
sheaves}'', {\bf Nagoya Math. Jour. 81}, (1981), 153-175.

\bibitem[Muk92]{muk92}
S. Mukai, 
''{\em Fano 3-folds}'', 

\bibitem[MumLC]{mumlc}
D. Mumford, 
''{\em Lectures on curves  on an algebraic surface}''. 
{\bf Annals of Math. Studies, n. 59, Princeton Univ. Press}, (1966).


\bibitem[Mum66-7]{mum66-7}
D. Mumford, 
''{\em On the equations defining abelian varieties.I, II, III}'',  
{ \bf Invent. Math. }1 (1966 )287--354,  3 (1967) , 75--135  , 3 (1967), 215--244. 



\bibitem[Mum68]{mum68}
D. Mumford, 
''{\em Rational equivalence of $0$-cycles on surfaces}''. 
{\bf J. Math. Kyoto Univ. 9}, (1968) 195--204.



\bibitem[Mum70]{mum70}
D. Mumford, 
''{\em Abelian varieties}''
{\bf Tata Institute of Fundamental Research Studies in Mathematics, 
 Oxford University Press. VIII} (1970).

\bibitem[Mum75]{mum74}
D. Mumford, 
''{\em Curves and their Jacobians}''
{\bf The University of Michigan Press, Ann. Arbor}, (1975) 104 p. . 
 
 \bibitem[Mum3]{mum3}
D. Mumford, 
''{\em Tata lectures on theta. I}'' With the assistance of C. Musili, M.
Nori, E. Previato and M. Stillman.
   {\bf Progress in Mathematics, 28. BirkhŠuser Boston, Inc., Boston, MA}, (1983).
xiii+235 pp.

 ''{\em Tata lectures on theta. II. 
   Jacobian theta functions and differential equations}''. With the collaboration of
C. Musili, M. Nori, E. Previato, M. Stillman and H. Umemura.{\bf Progress in
Mathematics, 43. BirkhŠuser Boston, Inc.,
   Boston, MA}, (1984) xiv+272 pp..
 
''{\em Tata lectures on theta. III}''
  With the collaboration of Madhav Nori and Peter Norman. {\bf Progress in
Mathematics, 97. BirkhŠuser Boston, Inc., Boston, MA} (1991), viii+202 pp.


s

\bibitem [N-D79]{n-d79}
 
 Norman, A.C.; Davenport, J.H.
``{\em Symbolic integration - the dust settles? }''
Symbolic and algebraic computation, EUROSAM '79, int. Symp., 
Marseille 1979, 
{\bf Lect. Notes Comput. Sci. 72}, 398-407 (1979).

\bibitem[Painl03]{painl}
P.Painlev\'e,
``{\em  Sur les fonctions qui admettent un th\'eoreme d' addition}'',
{\bf Acta Math., 27 }, (1903).

\bibitem[Pat99]{pat99}
 Patterson, S.J.
``{\em Abel's theorem. }''
Bottazzini, U. (ed.), Studies in the history of modern mathematics. 
IV. Palermo: {\bf Circolo Matemˆtico di Palermo, Suppl. Rend. Circ. Mat. Palermo, II.
Ser. 61}, 9-48 (1999)

\bibitem[Pic95]{pic}
\'E. Picard,
``{\em Sur la th\'eorie des groupes et des surfaces alg\'ebriques. }''
{\bf Rendiconti del Circ. Mat. di Palermo, IX } (1895).

\bibitem[Pir01]{pi01}
P. Pirola,
``{\em Algebraic surfaces with $p_g-q-3$ and no irrational pencils }'',
{\bf  Manuscripta Math. 108 } (2002), n.2, 163-170.


\bibitem[Poi84]{poi84}
H. Poincar\'e,
``{\em Sur la r\'eduction des int\'egrales ab\'eliennes. }''
{\bf Bull. Soc. Math. France 12 } (1884), 124-143.

\bibitem[Poi02]{poi02}
H. Poincar\'e,
``{\em Sur les fonctions ab\'eliennes. }''
{\bf Acta Math. 26 } (1902), 43- 98.

\bibitem[Poi10]{poi10}
H. Poincar\'e,
``{\em Sur les courbes trac\'ees sur les surfaces algebriques. }''
{\bf Annales de l' \'Ecole normale Sup.,III, 27 } (1910), 43- 98.



\bibitem[Ran81]{ran}
 R.H. Risch,  ``{\em On subvarieties of abelian varieties.}''
{\bf Inv. Math. 62},  (1981), 459-479.

\bibitem[RA]{RA}
 U. Bottazzini, A. Conte, P. Gario editors, F. Enriques , author,  ``{\em Riposte
Armonie.}'', {\bf Bollati Boringhieri, Torino}  (1996), 722 pp. .

\bibitem[Rec74]{rec}
S. Recillas, 
``{\em Jacobians of curves with $g^1_4 $'s are the Prym's of trigonal curves.}''
{\bf Bol. Soc. Mat. Mex., II. Ser. 19}, (1974)  9-13 .

\bibitem[Risch]{risch}
 R.H. Risch,  ``{\em The solution of the problem of integration in finite terms.}''
{\bf Bull. Amer. Math. Soc. 76},  (1970), 605--608.

\bibitem[Ros29]{ros}
C. Rosati,
``{\em Sulle matrici di Riemann. }''
{\bf  Rend. Circ. Mat. Palermo, 53 }, (1929), 79-134.

\bibitem[Schap90]{schap90}
N. Schappacher, 
``{\em D\'eveloppement de la loi de groupe sur une cubique. [Development of
the law of composition on a cubic curve] }''  {\bf S\'eminaire de ThŽorie des
Nombres, Paris 1988--1989,   Progr. Math., 91, 
BirkhŠuser Boston, Boston, MA},( 1990), 159--184. 

\bibitem[Scha-Scho]{scha-scho}
N.  Schappacher, R. Schoof,``{\em Beppo Levi and the arithmetic of elliptic
curves. }'' {\bf Math. Intelligencer 18 }(1996), no. 1, 57--69.


\bibitem[Schr]{schr}
F.O.Schreyer,
``{\em  Geometry and algebra of prime Fano 3-folds of genus $22$}'',
{\bf Comp. Math. }, (2002)  .


\bibitem[Scor99]{scor99}
G. Scorza, 
``{\em  Sopra la teoria delle figure polari delle curve piane del $4^o$
ordine}'', {\bf Ann.di Mat.(3) 2 }, (1899), 155-202.

\bibitem[Scor16]{scor16}
G. Scorza, 
``{\em  Intorno alla teoria generale delle matrici di Riemann}'', {\bf Rend. Circ.
Mat. Palermo 41 }, (1916), 263-379.

\bibitem[Sev03]{sev03}
F. Severi, 
``{\em Sulle superficie che rappresentano le coppie di punti  di una
curva algebrica.}''  {\bf Atti della R. Accad. delle Scienze di Torino, 38},
(1903), 185-200.



\bibitem[Sev05]{sev05}
F. Severi, 
``{\em Il teorema d' Abel sulle  superficie algebriche.}''  {\bf Annali di Mat.
s.III, t.XII}, (1905), 55-79.



\bibitem[Sev32]{sev32}
F. Severi, 
``{\em La serie canonica e la teoria delle serie principali di gruppi di punti sopra
una superficie algebrica.}''  {\bf Comment. Math. Helv. 4},  (1932), 268-326.

\bibitem[Sev47]{sev47}
F. Severi, 
``{\em 
Funzioni quasi-abeliane. }'' 
{\bf Pontificiae Academiae Scientiarum Scripta Varia 4. 
Roma: Pontificia Academia Scientiarum. }327 p. (1947). 

\bibitem[Sev-SEL]{sev-sel}
F. Severi, 
``{\em 
Memorie scelte. Vol. I.  }'' 
{\bf  Cesare Zuffi Editore. XX, Bologna}, 458 p. (1950). 

\bibitem[Shim63]{shim}
G. Shimura, ''{\em On analytic families of polarized abelian varieties and
automorphic functions}'', {\bf Ann. of Math. 78 }, (1963), 149-193.

\bibitem[Sieg43]{sieg43}
C.L. Siegel, ''{\em Symplectic geometry}'' 
{\bf Amer. Jour. Math., n. 65} (1943), 1-86, reprinted by Academic Press, New
York and London, (1964).

\bibitem[Sieg73]{sieg}
C.L. Siegel, ''{\em Topics in complex function theory, Vol. III}'' 
{\bf Interscience Tracts in pure and appl. Math., n. 25, Wiley, N.Y.} (1973).

\bibitem[Siu02]{siu02}
J. T. Siu, 
''{\em  Extension of twisted pluricanonical sections with plurisubharmonic weight
and invariance of semipositively twisted plurigenera for manifolds not
necessarily of general type. }'' in 'Complex Geometry, Collections of papers
dedicated to Hans Grauert, I. Bauer et al ed. {\bf Springer Verlag}, 223-277
(2002).

\bibitem[Tate74]{tate}
J. T. Tate, 
''{\em  The arithmetic of elliptic curves. }'' 
{\bf Invent. Math. 23}, 179-206 (1974).

\bibitem[Tor13]{tor}
R. Torelli, 
''{\em  Sulle varieta' di Jacobi, I, II. }'' 
{\bf Rendiconti R. Accad. dei Lincei 22-2}, (1913), 98-103, 437-441.


\bibitem[Trag79]{trag79}
B.M. Trager, 
''{\em Integration of simple radical extensions. }'' 
Symbolic and algebraic computation (EUROSAM '79, Internat. Sympos., Marseille, 1979), pp. 408--414, 
{\bf Lecture Notes in Comput. Sci., 72}, 
Springer, Berlin-New York, (1979). 

\bibitem[Tric]{tric}
F. Tricomi, 
''{\em Funzioni ellittiche.  }'' 
Second edition, 
{\bf Nicola Zanichelli Editore, Bologna}, (1951), ix+343 pp.

\bibitem [Ue75]{ue75} 
  K. Ueno, ``{\em
Classification theory of algebraic varieties and compact complex spaces.}'' 

{\bf Lecture Notes in Mathematics 439,} 
Springer-Verlag, XIX, 278 p.(1975).

\bibitem [Weier78]{weier} 

K. Weierstra\ss , 
``{\em Einleitung in die Theorie der analytischen Funktionen. }'' , Vorlesung Berlin
1878, Lecture notes taken and with supplementary material
   by Adolf Hurwitz. Preface by R. Remmert, ed.  Peter
Ullrich. {\bf Dokumente zur Geschichte der Mathematik  4. Deutsche Mathematiker
Vereinigung, Freiburg; Friedr. Vieweg
  Sohn, Braunschweig}, (1988) xxx+184 pp..

\bibitem[Weil-29]{weil29}
A. Weil, 
``{\em 
L' arithm\'etique sur les courbes alg\'ebriques. }'' 
{ \bf Acta Math. 52}, (1929), 281-315.

\bibitem[Weil-CA]{weilca}
A. Weil, 
``{\em 
Sur les courbes alg\'ebriques et les variŽtŽs qui s'en deduisent. }'' 
{\bf Actualit\'es scientifiques et industrielles. 1041. Paris: Hermann \& Cie.
IV}, 85 p. (1948); Publ. Inst. Math. Univ. Strasbourg 7 (1945). 

\bibitem[Weil-VA]{weilva}
A. Weil, 
``{\em Vari\'et\'es abeliennes et courbes alg\'ebriques. }'' 
{\bf Paris: Hermann \& Cie.} 163 p. (1948).

\bibitem[Weil-57]{weil57}
A. Weil, 
``{\em 
Zum Beweis des Torellischen Satzes. }'' 
{ \bf Nachr. Akad. Wiss. Gšttingen, Math.-Phys. Kl., IIa 1957}, 32-53 (1957).

\bibitem[Weil-64]{weil64}
A. Weil, 
``{\em 
Sur certains groupes d'op\'erateurs unitaires. }''
{\bf Acta Math. 111}, (1964), 143-211 . 

\bibitem[Weil-VK]{weilVK}
A. Weil, 
``{\em Introduction a l' \'etude des vari\'et\'es k\"ahl\'eriennes. }''  {\bf
Paris: Hermann \& Cie.}, (1958).


\bibitem[Weyl31]{Weyl31}

H.Weyl, ``{\em  Topologie und abstrakte Algebra als zwei Wege mathematischen
Verst\"andnisses}'' {\bf  Unterrichtsbl\"atter f\"ur Mathematik und
Naturwissenschaften, band 38 }, (1933), s. 177-188, also in 

\bibitem[Weyl]{Weyl}
H.Weyl, ``{\em
Gesammelte Abhandlungen. B\"ande I, II, III, IV.}'' 
Edited by K. Chandrasekharan 
{\bf Springer-Verlag, Berlin-New York } (1968) ,Band I: vi+698 pp.,
Band II: iv+647 pp., Band III: iv+791 pp., Band IV: ix+694 pp..

\bibitem[Weyl34]{Weyl34}

H.Weyl, ``{\em  On generalized Riemann matrices}'', {\bf Ann.Math. 35 },
(1934), 714-729.

\bibitem[Weyl36]{Weyl36}

H.Weyl, ``{\em  Generalized Riemann matrices and factor sets}'', {\bf Ann.Math. 37
}, (1936), 709-745.


\bibitem[Yag]{yag}
 I. M.  Yaglom,``{\em  Felix Klein and Sophus Lie. Evolution of the idea of
symmetry in the nineteenth century.}'' Translated from the Russian by Sergei
Sossinsky. Translation edited by Hardy Grant and Abe Shenitzer.{\bf BirkhŠuser
Boston, Inc.,
   Boston, MA} (1988). xii+237 pp..

\bibitem[Zap91]{zap}
  G.  Zappa, 
``{\em The papers of Gaetano Scorza on group theory. }''
{ \bf Atti Accad. Naz. Lincei,
Cl. Sci. Fis. Mat. Nat., IX. Ser., Rend. Lincei, Mat. Appl. 2, No.2}, (1991), 95-101 .

\bibitem[Zar35]{zar}
O. Zariski,
``{\em  Algebraic surfaces}'', {\bf Ergebnisse der Math. und ihrer Grenzgebiete,
Band III Heft 5  , Springer Verlag} (1935), second supplemented edition in {\bf
Ergebnisse 61}, (1970), 270 pp.






\end{thebibliography}
\end{document}